
\def\Resetstrings{
    \def\present{ }\let\bgroup={\let\egroup=}
    \def\Astr{}\def\astr{}\def\Atest{}\def\atest{}%
    \def\Bstr{}\def\bstr{}\def\Btest{}\def\btest{}%
    \def\Cstr{}\def\cstr{}\def\Ctest{}\def\ctest{}%
    \def\Dstr{}\def\dstr{}\def\Dtest{}\def\dtest{}%
    \def\Estr{}\def\estr{}\def\Etest{}\def\etest{}%
    \def\Fstr{}\def\fstr{}\def\Ftest{}\def\ftest{}%
    \def\Gstr{}\def\gstr{}\def\Gtest{}\def\gtest{}%
    \def\Hstr{}\def\hstr{}\def\Htest{}\def\htest{}%
    \def\Istr{}\def\istr{}\def\Itest{}\def\itest{}%
    \def\Jstr{}\def\jstr{}\def\Jtest{}\def\jtest{}%
    \def\Kstr{}\def\kstr{}\def\Ktest{}\def\ktest{}%
    \def\Lstr{}\def\lstr{}\def\Ltest{}\def\ltest{}%
    \def\Mstr{}\def\mstr{}\def\Mtest{}\def\mtest{}%
    \def\Nstr{}\def\nstr{}\def\Ntest{}\def\ntest{}%
    \def\Ostr{}\def\ostr{}\def\Otest{}\def\otest{}%
    \def\Pstr{}\def\pstr{}\def\Ptest{}\def\ptest{}%
    \def\Qstr{}\def\qstr{}\def\Qtest{}\def\qtest{}%
    \def\Rstr{}\def\rstr{}\def\Rtest{}\def\rtest{}%
    \def\Sstr{}\def\sstr{}\def\Stest{}\def\stest{}%
    \def\Tstr{}\def\tstr{}\def\Ttest{}\def\ttest{}%
    \def\Ustr{}\def\ustr{}\def\Utest{}\def\utest{}%
    \def\Vstr{}\def\vstr{}\def\Vtest{}\def\vtest{}%
    \def\Wstr{}\def\wstr{}\def\Wtest{}\def\wtest{}%
    \def\Xstr{}\def\xstr{}\def\Xtest{}\def\xtest{}%
    \def\Ystr{}\def\ystr{}\def\Ytest{}\def\ytest{}%
}
\Resetstrings

\def\Refformat{
         \if\Jtest\present
             {\if\Vtest\present\journalarticleformat
                  \else\conferencereportformat\fi}
            \else\if\Btest\present\bookarticleformat
               \else\if\Rtest\present\technicalreportformat
                  \else\if\Itest\present\bookformat
                     \else\otherformat\fi\fi\fi\fi}

\def\Rpunct{
   \def\Lspace{ }%
   \def\Lperiod{ }
   \def\Lcomma{ }
   \def\Lquest{ }
   \def\Lcolon{ }
   \def\Lscolon{ }
   \def\Lbang{ }
   \def\Lquote{ }
   \def\Lqquote{ }
   \def\Lrquote{ }
   \def\Rspace{}%
   \def\Rperiod{.}
   \def\Rcomma{,}
   \def\Rquest{?}
   \def\Rcolon{:}
   \def\Rscolon{;}
   \def\Rbang{!}
   \def\Rquote{'}
   \def\Rqquote{"}
   \def\Rrquote{`}
   }

\def\Lpunct{
   \def\Lspace{}%
   \def\Lperiod{\unskip.}
   \def\Lcomma{\unskip,}
   \def\Lquest{\unskip?}
   \def\Lcolon{\unskip:}
   \def\Lscolon{\unskip;}
   \def\Lbang{\unskip!}
   \def\Lquote{\unskip'}
   \def\Lqquote{\unskip"}
   \def\Lrquote{\unskip`}
   \def\Rspace{\spacefactor=1000}%
   \def\Rperiod{\spacefactor=3000}
   \def\Rcomma{\spacefactor=1250}
   \def\Rquest{\spacefactor=3000}
   \def\Rcolon{\spacefactor=2000}
   \def\Rscolon{\spacefactor=1250}
   \def\Rbang{\spacefactor=3000}
   \def\Rquote{\spacefactor=1000}
   \def\Rqquote{\spacefactor=1000}
   \def\Rrquote{\spacefactor=1000}
   }

\def\Refstd{
     \def\Acomma{\unskip, }
     \def\Aand{\unskip\ and }
     \def\Aandd{\unskip\ and }
     \def\Ecomma{\unskip, }
     \def\Eand{\unskip\ and }
     \def\Eandd{\unskip\ and }
     \def\acomma{\unskip, }
     \def\aand{\unskip\ and }
     \def\aandd{\unskip\ and }
     \def\ecomma{\unskip, }
     \def\eand{\unskip\ and }
     \def\eandd{\unskip\ and }
     \def\Namecomma{\unskip, }
     \def\Nameand{\unskip\ and }
     \def\Nameandd{\unskip\ and }
     \def\Revcomma{\unskip, }
     \def\Initper{.\ }
     \def\Initgap{\dimen0=\spaceskip\divide\dimen0 by 2\hskip-\dimen0}%
   }

\def\Smallcapsaand{
     \def\Aand{\unskip\bgroup{\Smallcapsfont\ AND }\egroup}%
     \def\Aandd{\unskip\bgroup{\Smallcapsfont\ AND }\egroup}%
     \def\eand{\unskip\bgroup\Smallcapsfont\ AND \egroup}%
     \def\eandd{\unskip\bgroup\Smallcapsfont\ AND \egroup}%
   }

\def\Smallcapseand{
     \def\Eand{\unskip\bgroup\Smallcapsfont\ AND \egroup}%
     \def\Eandd{\unskip\bgroup\Smallcapsfont\ AND \egroup}%
     \def\aand{\unskip\bgroup\Smallcapsfont\ AND \egroup}%
     \def\aandd{\unskip\bgroup\Smallcapsfont\ AND \egroup}%
   }

   \def\Citefont{}
   \def\ACitefont{}
   \def\Authfont{}
   \def\Titlefont{}
   \def\Tomefont{\sl}
   \def\Volfont{}
   \def\Flagfont{}
   \def\Reffont{\rm}
   \def\Smallcapsfont{\sevenrm}
   \def\Flagstyle#1{\hangindent\parindent\indent\hbox to0pt
       {\hss[{\Flagfont#1}]\kern.5em}\ignorespaces}


\def\Citebrackets{\Rpunct
   \def\Lcitemark{\def\Cfont{\Citefont}[\bgroup\Cfont}
   \def\Rcitemark{\egroup]}
   \def\LAcitemark{\def\Cfont{\ACitefont}\bgroup\ACitefont}%
   \def\RAcitemark{\egroup}
   \def\LIcitemark{\egroup}
   \def\RIcitemark{\bgroup\Cfont}
   \def\Citehyphen{\egroup--\bgroup\Cfont}
   \def\Citecomma{\egroup,\hskip0pt\bgroup\Cfont}%
   \def\Citebreak{}
   }

\def\Citeparen{\Rpunct
   \def\Lcitemark{\def\Cfont{\Citefont}(\bgroup\Cfont}
   \def\Rcitemark{\egroup)}
   \def\LAcitemark{\def\Cfont{\ACitefont}\bgroup\ACitefont}%
   \def\RAcitemark{\egroup}
   \def\LIcitemark{\egroup}
   \def\RIcitemark{\bgroup\Cfont}
   \def\Citehyphen{\egroup--\bgroup\Cfont}
   \def\Citecomma{\egroup,\hskip0pt\bgroup\Cfont}%
   \def\Citebreak{}
   }

\def\Citesuper{\Lpunct
   \def\Lcitemark{\def\Cfont{\Citefont}\raise1ex\hbox\bgroup\bgroup\Cfont}%
   \def\Rcitemark{\egroup\egroup}
   \def\LAcitemark{\def\Cfont{\ACitefont}\bgroup\ACitefont}%
   \def\RAcitemark{\egroup}
   \def\LIcitemark{\egroup\egroup}
   \def\RIcitemark{\raise1ex\hbox\bgroup\bgroup\Cfont}%
   \def\Citehyphen{\egroup--\bgroup\Cfont}
   \def\Citecomma{\egroup,\hskip0pt\bgroup%
      \Cfont}
   \def\Citebreak{}
   }

\def\Citenamedate{\Rpunct
   \def\Lcitemark{
      \def\Citebreak{\egroup\ [\bgroup\Citefont}
      \def\Citecomma{\egroup]; 
         \bgroup\let\uchyph=1\Citefont}(\bgroup\let\uchyph=1\Citefont}%
   \def\Rcitemark{\egroup])}
   \def\LAcitemark{
      \def\Citebreak{\egroup\ [\bgroup\Citefont}\def\Citecomma{\egroup], %
         \bgroup\ACitefont }\bgroup\let\uchyph=1\ACitefont}%
   \def\RAcitemark{\egroup]}
  \def\Citehyphen{\egroup--\bgroup\Citefont}
   \def\LIcitemark{\egroup}
   \def\RIcitemark{\bgroup\Citefont}
   }

\def\Flagstyle#1{\Flagfont#1. }

\Refstd\Citebrackets
\def\Citefont{\bf}\def\Titlefont{\sl}\def\Volfont{\bf}\def\Tomefont{\Reffont}

\def\journalarticleformat{\Reffont\let\uchyph=1\parindent=1.25pc\def\Comma{}%

\sfcode`\.=1000\sfcode`\?=1000\sfcode`\!=1000\sfcode`\:=1000\sfcode`\;=1000\sfcode`\,=1000
                \par\vfil\penalty-200\vfilneg
      \if\Ftest\present\Flagstyle\Fstr\fi%
       \if\Atest\present\bgroup\Authfont\Astr\egroup\def\Comma{\unskip, }\fi%
        \if\Ttest\present\Comma\bgroup\Titlefont\Tstr\egroup\def\Comma{, }\fi%
         \if\etest\present\hskip.2em(\bgroup\estr\egroup)\def\Comma{\unskip,
}\fi%
          \if\Jtest\present\Comma\bgroup\Tomefont\Jstr\/\egroup\def\Comma{,
}\fi%

\if\Vtest\present\if\Jtest\present\hskip.2em\else\Comma\fi\bgroup\Volfont\Vstr\egroup\def\Comma{, }\fi%
            \if\Dtest\present\hskip.2em(\bgroup\Dstr\egroup)\def\Comma{, }\fi%
             \if\Ptest\present\bgroup, \Pstr\egroup\def\Comma{, }\fi%
              \if\ttest\present\Comma\bgroup\Titlefont\tstr\egroup\def\Comma{,
}\fi%

\if\jtest\present\Comma\bgroup\Tomefont\jstr\/\egroup\def\Comma{, }\fi%

\if\vtest,\present\if\jtest\present\hskip.2em\else\Comma\fi\bgroup\Volfont\vstr\egroup\def\Comma{, }\fi%
                 \if\dtest\present\hskip.2em(\bgroup\dstr\egroup)\def\Comma{,
}\fi%
                  \if\ptest\present\bgroup, \pstr\egroup\def\Comma{, }\fi%
                   \if\Gtest\present{\Comma Gov't ordering no.
}\bgroup\Gstr\egroup\def\Comma{, }\fi%
                    \if\Mtest\present\Comma MR
\#\bgroup\Mstr\egroup\def\Comma{, }\fi%
                     \if\Otest\present{\Comma\bgroup\Ostr\egroup.}\else{.}\fi%
                      \vskip3ptplus1ptminus1pt}

\def\conferencereportformat{\Reffont\let\uchyph=1\parindent=1.25pc\def\Comma{}%

\sfcode`\.=1000\sfcode`\?=1000\sfcode`\!=1000\sfcode`\:=1000\sfcode`\;=1000\sfcode`\,=1000
                \par\vfil\penalty-200\vfilneg
      \if\Ftest\present\Flagstyle\Fstr\fi%
       \if\Atest\present\bgroup\Authfont\Astr\egroup\def\Comma{\unskip, }\fi%
        \if\Ttest\present\Comma\bgroup\Titlefont\Tstr\egroup\def\Comma{, }\fi%
         \if\Jtest\present\Comma\bgroup\Tomefont\Jstr\/\egroup\def\Comma{,
}\fi%
          \if\Ctest\present\Comma\bgroup\Cstr\egroup\def\Comma{, }\fi%
           \if\Dtest\present\hskip.2em(\bgroup\Dstr\egroup)\def\Comma{, }\fi%
            \if\Mtest\present\Comma MR \#\bgroup\Mstr\egroup\def\Comma{, }\fi%
             \if\Otest\present{\Comma\bgroup\Ostr\egroup.}\else{.}\fi%
              \vskip3ptplus1ptminus1pt}

\def\bookarticleformat{\Reffont\let\uchyph=1\parindent=1.25pc\def\Comma{}%

\sfcode`\.=1000\sfcode`\?=1000\sfcode`\!=1000\sfcode`\:=1000\sfcode`\;=1000\sfcode`\,=1000
                \par\vfil\penalty-200\vfilneg
      \if\Ftest\present\Flagstyle\Fstr\fi%
       \if\Atest\present\bgroup\Authfont\Astr\egroup\def\Comma{\unskip, }\fi%
        \if\Ttest\present\Comma\bgroup\Titlefont\Tstr\egroup\def\Comma{, }\fi%
         \if\etest\present\hskip.2em(\bgroup\estr\egroup)\def\Comma{\unskip,
}\fi%
          \if\Btest\present\Comma in
\bgroup\Tomefont\Bstr\/\egroup\def\Comma{\unskip, }\fi%
           \if\otest\present\ \bgroup\ostr\egroup\def\Comma{, }\fi%
            \if\Etest\present\Comma\bgroup\Estr\egroup\unskip,
\ifnum\Ecnt>1eds.\else ed.\fi\def\Comma{, }\fi%
             \if\Stest\present\Comma\bgroup\Sstr\egroup\def\Comma{, }\fi%
              \if\Vtest\present\Comma vol. \bgroup\Vstr\egroup\def\Comma{,
}\fi%
               \if\Ntest\present\Comma no. \bgroup\Nstr\egroup\def\Comma{,
}\fi%
                \if\Itest\present\Comma\bgroup\Istr\egroup\def\Comma{, }\fi%
                 \if\Ctest\present\Comma\bgroup\Cstr\egroup\def\Comma{, }\fi%
                  \if\Dtest\present\Comma\bgroup\Dstr\egroup\def\Comma{, }\fi%
                   \if\Ptest\present\Comma\Pstr\def\Comma{, }\fi%

\if\ttest\present\Comma\bgroup\Titlefont\Tstr\egroup\def\Comma{, }\fi%
                     \if\btest\present\Comma in
\bgroup\Tomefont\bstr\egroup\def\Comma{, }\fi%
                       \if\atest\present\Comma\bgroup\astr\egroup\unskip,
\if\acnt\present eds.\else ed.\fi\def\Comma{, }\fi%
                        \if\stest\present\Comma\bgroup\sstr\egroup\def\Comma{,
}\fi%
                         \if\vtest\present\Comma vol.
\bgroup\vstr\egroup\def\Comma{, }\fi%
                          \if\ntest\present\Comma no.
\bgroup\nstr\egroup\def\Comma{, }\fi%

\if\itest\present\Comma\bgroup\istr\egroup\def\Comma{, }\fi%

\if\ctest\present\Comma\bgroup\cstr\egroup\def\Comma{, }\fi%

\if\dtest\present\Comma\bgroup\dstr\egroup\def\Comma{, }\fi%
                              \if\ptest\present\Comma\pstr\def\Comma{, }\fi%
                               \if\Gtest\present{\Comma Gov't ordering no.
}\bgroup\Gstr\egroup\def\Comma{, }\fi%
                                \if\Mtest\present\Comma MR
\#\bgroup\Mstr\egroup\def\Comma{, }\fi%

\if\Otest\present{\Comma\bgroup\Ostr\egroup.}\else{.}\fi%
                                  \vskip3ptplus1ptminus1pt}

\def\bookformat{\Reffont\let\uchyph=1\parindent=1.25pc\def\Comma{}%

\sfcode`\.=1000\sfcode`\?=1000\sfcode`\!=1000\sfcode`\:=1000\sfcode`\;=1000\sfcode`\,=1000
                \par\vfil\penalty-200\vfilneg
      \if\Ftest\present\Flagstyle\Fstr\fi%
       \if\Atest\present\bgroup\Authfont\Astr\egroup\def\Comma{\unskip, }%

\else\if\Etest\present\bgroup\def\Eand{\Aand}\def\Eandd{\Aandd}\Authfont\Estr\egroup\unskip, \ifnum\Ecnt>1eds.\else ed.\fi\def\Comma{, }%

\else\if\Itest\present\bgroup\Authfont\Istr\egroup\def\Comma{, }\fi\fi\fi%

\if\Ttest\present\Comma\bgroup\Titlefont\Tstr\/\egroup\def\Comma{\unskip, }%

\else\if\Btest\present\Comma\bgroup\Titlefont\Bstr\/\egroup\def\Comma{\unskip,
}\fi\fi%
            \if\otest\present\ \bgroup\ostr\egroup\def\Comma{, }\fi%

\if\etest\present\hskip.2em(\bgroup\estr\egroup)\def\Comma{\unskip, }\fi%
              \if\Stest\present\Comma\bgroup\Sstr\egroup\def\Comma{, }\fi%
               \if\Vtest\present\Comma vol. \bgroup\Vstr\egroup\def\Comma{,
}\fi%
                \if\Ntest\present\Comma no. \bgroup\Nstr\egroup\def\Comma{,
}\fi%
                 \if\Atest\present\if\Itest\present
                         \Comma\bgroup\Istr\egroup\def\Comma{\unskip, }\fi%
                      \else\if\Etest\present\if\Itest\present
                              \Comma\bgroup\Istr\egroup\def\Comma{\unskip,
}\fi\fi\fi%
                     \if\Ctest\present\Comma\bgroup\Cstr\egroup\def\Comma{,
}\fi%
                      \if\Dtest\present\Comma\bgroup\Dstr\egroup\def\Comma{,
}\fi%

\if\ttest\present\Comma\bgroup\Titlefont\tstr\egroup\def\Comma{, }%

\else\if\btest\present\Comma\bgroup\Titlefont\bstr\egroup\def\Comma{, }\fi\fi%

\if\stest\present\Comma\bgroup\sstr\egroup\def\Comma{, }\fi%
                           \if\vtest\present\Comma vol.
\bgroup\vstr\egroup\def\Comma{, }\fi%
                            \if\ntest\present\Comma no.
\bgroup\nstr\egroup\def\Comma{, }\fi%

\if\itest\present\Comma\bgroup\istr\egroup\def\Comma{, }\fi%

\if\ctest\present\Comma\bgroup\cstr\egroup\def\Comma{, }\fi%

\if\dtest\present\Comma\bgroup\dstr\egroup\def\Comma{, }\fi%
                                \if\Gtest\present{\Comma Gov't ordering no.
}\bgroup\Gstr\egroup\def\Comma{, }\fi%
                                 \if\Mtest\present\Comma MR
\#\bgroup\Mstr\egroup\def\Comma{, }\fi%

\if\Otest\present{\Comma\bgroup\Ostr\egroup.}\else{.}\fi%
                                   \vskip3ptplus1ptminus1pt}

\def\technicalreportformat{\Reffont\let\uchyph=1\parindent=1.25pc\def\Comma{}%

\sfcode`\.=1000\sfcode`\?=1000\sfcode`\!=1000\sfcode`\:=1000\sfcode`\;=1000\sfcode`\,=1000
                \par\vfil\penalty-200\vfilneg
      \if\Ftest\present\Flagstyle\Fstr\fi%
       \if\Atest\present\bgroup\Authfont\Astr\egroup\def\Comma{\unskip, }%

\else\if\Etest\present\bgroup\def\Eand{\Aand}\def\Eandd{\Aandd}\Authfont\Estr\egroup\unskip, \ifnum\Ecnt>1eds.\else ed.\fi\def\Comma{, }%

\else\if\Itest\present\bgroup\Authfont\Istr\egroup\def\Comma{, }\fi\fi\fi%
          \if\Ttest\present\Comma\bgroup\Titlefont\Tstr\egroup\def\Comma{,
}\fi%
           \if\Atest\present\if\Itest\present
                   \Comma\bgroup\Istr\egroup\def\Comma{, }\fi%
                \else\if\Etest\present\if\Itest\present
                        \Comma\bgroup\Istr\egroup\def\Comma{, }\fi\fi\fi%
            \if\Rtest\present\Comma\bgroup\Rstr\egroup\def\Comma{, }\fi%
             \if\Ctest\present\Comma\bgroup\Cstr\egroup\def\Comma{, }\fi%
              \if\Dtest\present\Comma\bgroup\Dstr\egroup\def\Comma{, }\fi%
               \if\ttest\present\Comma\bgroup\Titlefont\tstr\egroup\def\Comma{,
}\fi%
                \if\itest\present\Comma\bgroup\istr\egroup\def\Comma{, }\fi%
                 \if\rtest\present\Comma\bgroup\rstr\egroup\def\Comma{, }\fi%
                  \if\ctest\present\Comma\bgroup\cstr\egroup\def\Comma{, }\fi%
                   \if\dtest\present\Comma\bgroup\dstr\egroup\def\Comma{, }\fi%
                    \if\Gtest\present{\Comma Gov't ordering no.
}\bgroup\Gstr\egroup\def\Comma{, }\fi%
                     \if\Mtest\present\Comma MR
\#\bgroup\Mstr\egroup\def\Comma{, }\fi%
                      \if\Otest\present{\Comma\bgroup\Ostr\egroup.}\else{.}\fi%
                       \vskip3ptplus1ptminus1pt}

\def\otherformat{\Reffont\let\uchyph=1\parindent=1.25pc\def\Comma{}%

\sfcode`\.=1000\sfcode`\?=1000\sfcode`\!=1000\sfcode`\:=1000\sfcode`\;=1000\sfcode`\,=1000
                \par\vfil\penalty-200\vfilneg
      \if\Ftest\present\Flagstyle\Fstr\fi%
       \if\Atest\present\bgroup\Authfont\Astr\egroup\def\Comma{\unskip, }%

\else\if\Etest\present\bgroup\def\Eand{\Aand}\def\Eandd{\Aandd}\Authfont\Estr\egroup\unskip, \ifnum\Ecnt>1eds.\else ed.\fi\def\Comma{, }%

\else\if\Itest\present\bgroup\Authfont\Istr\egroup\def\Comma{, }\fi\fi\fi%
          \if\Ttest\present\Comma\bgroup\Titlefont\Tstr\egroup\def\Comma{,
}\fi%
            \if\Atest\present\if\Itest\present
                    \Comma\bgroup\Istr\egroup\def\Comma{, }\fi%
                 \else\if\Etest\present\if\Itest\present
                         \Comma\bgroup\Istr\egroup\def\Comma{, }\fi\fi\fi%
                 \if\Ctest\present\Comma\bgroup\Cstr\egroup\def\Comma{, }\fi%
                  \if\Dtest\present\Comma\bgroup\Dstr\egroup\def\Comma{, }\fi%
                   \if\Gtest\present{\Comma Gov't ordering no.
}\bgroup\Gstr\egroup\def\Comma{, }\fi%
                    \if\Mtest\present\Comma MR
\#\bgroup\Mstr\egroup\def\Comma{, }\fi%
                     \if\Otest\present{\Comma\bgroup\Ostr\egroup.}\else{.}\fi%
                      \vskip3ptplus1ptminus1pt}

\message {)}\message {DOCUMENT TEXT}
\magnification=\magstep1
\input amstex
\documentstyle{amsppt}
\voffset=-2cm
\topmatter
\title
Automorphic $L$\snug-functions, intertwining operators, and the irreducible
tempered representations of  $p$\snug-adic groups
\endtitle
\author
David Goldberg\footnotemark"*" and Freydoon Shahidi\footnotemark"**"\endauthor
\address
Department of Mathematics,
Purdue University,
West Lafayette, IN 47907
\endaddress
\leftheadtext{David Goldberg and Freydoon Shahidi}
\rightheadtext{$L$\snug-functions and intertwining operators}
\endtopmatter
\footnotetext""{*Partially supported by National Science Foundation
Postdoctoral Fellowship DMS9206246}
\footnotetext""{**Partially supported by National Science Foundation Grant
DMS9301040}

\def\BA{\Bbb A}

\def\calA{\Cal{A}}

\define\PPP{\text{\bf P}}

\def\G{\text{\bf G}}

\def\ds{\displaystyle}
\def\bC{\Bbb C}
\def\R{\Bbb R}
\def\bz{\Bbb Z}
\define\g{\frak g}
\def\afa{\frak a}
\def\CalE{\Cal E}
\def\bc{\Bbb C}
\def\a{\alpha}
\def\b{\beta}
\redefine\d{\delta}
\def\e{\varepsilon}
\def\l{\lambda}
\define\p{\psi}
\def\ph{\varphi}
\def\t{\theta}
\def\s{\sigma}
\def\l{\lambda}
\def\DE{\Delta}
\redefine\P{\Phi}
\define\lra{\longrightarrow}
\define\Ind{\operatorname{Ind}}

\def\A{\text{\bf A}}
\def\B{\text{\bf B}}
\def\TT{\text{\bf T}}
\def\M{\text{\bf M}}
\def\N{\text{\bf N}}
\def\U{\text{\bf U}}

\define\wt{\widetilde}

\baselineskip=18pt
\parskip=6pt
\NoBlackBoxes
\document
\baselineskip=18pt
\parskip=6pt
\subheading{Introduction}
The problem of identifying the tempered, or admissible dual of a reductive
group consists of two problems.  Determine the discrete series representations
of Levi subgroups of the group, and decompose the resulting parabolically
induced representations.
Neither problem is resolved in any generality, and while\Lspace \Lcitemark
11\Rcitemark \Rspace{} discusses the first,
here we concentrate on the second and point out recent developments on
the subject.

For the tempered dual, we will mainly discuss our approach, which is
based on the theory of $R$\snug-groups,
a method with application to the theory of automorphic forms\Lspace \Lcitemark
1\Citecomma
2\Rcitemark \Rspace{}.
Roughly speaking,  $R$\snug-groups are finite groups whose duals
parameterize irreducible constituents of representations parabolically induced
from discrete series, i.e., the non-discrete tempered spectrum.
To determine the $R$\snug-group, one needs to determine the zeros of
the Plancherel measure, a measure supported on the tempered spectrum
whose restriction to the discrete part gives their formal degrees.
On the other hand, a conjecture of Langlands relates the Plancherel measures to
certain objects of arithmetic significance.  This has played an important role
in the recent progress.
Our goal is to describe this crucial relationship between
arithmetic and harmonic analysis.

There are several reasonable expositions of the Langlands
program.  For the questions that we are considering, one
should consult\Lspace \Lcitemark 22\Rcitemark \Rspace{}.  We also call
attention
to\Lspace \Lcitemark 19\Rcitemark \Rspace{} for its clarity, as well as its
annotated
bibliography.

\subheading{\S 1 $L$\snug-functions and the Langlands program}
We give a brief introduction to the Langlands program, with an eye toward our
goal of describing the classification of irreducible tempered representations
of reductive groups over local fields.   For much more comprehensive
introductions, and more motivation, one should see
\Lcitemark 7\Citecomma
18\Citecomma
19\Citecomma
22\Citecomma
69\Rcitemark \Rspace{}.
For our purposes, it is enough to remark that the motivation for the Langlands
program is local and global class field theory.  Namely, the remarkable fact
that the characters (i.e. one dimensional representations) of  the
multiplicative group of a local (or global) field $F$ are naturally
parameterized by characters of the Galois group of $\bar F/F.$   Viewing
$F^\times$ as $GL_1(F),$ one hopes to describe the representations of
$\G(F)$ through some kind of Galois representation theory.

\subheading{(a) $L$--groups}
Let $F$ be a local non--archimedean field, of characteristic zero
and let
$\Cal O$ be the ring of integers of $F.$
We denote by $\frak P,$  the prime ideal of $F,$ and set
$q=|\Cal O/\frak{P}|.$
Choose  a  non-trivial unramified character $\psi$ of $F.$
Let
$\G$ be  a connected, reductive, algebraic group, defined over $F$
\Lcitemark 57\Rcitemark \Rspace{}.
We will assume that
\roster
\item $\bold{G}$ is quasi--split, i.e.\ there exists a Borel
subgroup $\bold{B}=\bold T\bold{U}$ defined over $F.$

\item $\bold{G}$ splits over an unramified extension $L/F,$
i.e., $\text{\bf T}(L)\simeq (L^\times)^{\dim\text{\bf T}}.$
\endroster

We recall from\Lspace \Lcitemark 57\Rcitemark \Rspace{} that $\bold{G}$ is
given by a
based root datum: $$\Psi=(X^{*}(\bold{T}), \
\triangle, \ X_{*}(\bold{T}), \check\triangle),$$ where $X^{*}(\bold{T})$
is the group of rational characters of $\bold T, \triangle$ a choice of simple
roots,
$X_{*}(\bold{T})$ the cocharacters, and $\check\triangle$ the simple
dual roots.
Let $\check\Psi$ be the based root datum given by
$\check\Psi=(X_*(\bold T), \check\triangle, X^{*}(\bold T),\triangle).$
Let ${}^L\! G^{0}$ be the
complex group with root datum $\check\Psi$
(see\Lspace \Lcitemark 57\LIcitemark{}, 4.11\RIcitemark \Rcitemark \Rspace{}).
Let $\Gamma_{L/F}$ be the
Galois group of $L/F.$
Define ${}^L\! G= {}^L\! G^{0}\ltimes\Gamma_{L/F}.$

Here, $\Gamma_{L/F}$ acts on ${}^L\! G^{0}$ by its action on the root datum
$\check\Psi.$  Note that whenever $G$ is a split group over $F,$
then $\Gamma_{L/F}$ acts trivially on the root datum $\Psi,$
and hence, ${}^L\! G= {}^L\! G^{0}\times\Gamma_{L/F}.$

\remark{Remark} In general $\Gamma_{L/F}$ must be replaced by
the Weil group (at least) but, for the moment,
we choose this form for simplicity.
Later,we will use the Weil group\Lspace \Lcitemark 80\Rcitemark \Rspace{}.
\endremark

\example{Examples:}

(1) $\bold{G}=GL_{n}, \quad \Psi=\check\Psi,$
so ${}^L\! G^{0} = GL_{n}(\Bbb C),$ and thus, $^L\! G = GL_{n}(\Bbb
C)\times\Gamma_{L/F}.$  Note that if $n=1,$ then $^L\!G^0=\bc^\times.$
Thus, local class field theory says that $\hat G(F)$ is parameterized by the
homomorphisms from $\Gamma_{\bar F/F}$ to $^L\!G^0.$

(2) $\G = Sp_{2n} = \{g\in GL_{2n}|{}^t\! gJg= J\},$
where $J,$ a symplectic
form.
For the purpose of this example, we choose
$J=\pmatrix \hfill 0 & \hfill I_{n}\\
-I_{n} & \hfill 0\endpmatrix.$
$$
\bold{T}=\left\{
\pmatrix x_{1}&&&&&&&\\
& x_{2}&&&&&&\\
& & \ddots & & 0\\
&&& x_{n}&&&&\\
& 0&&&x^{-1}_{1}&&&\\
&&&&&\ddots&&&\\
&&&&&&& x^{-1}_{n}\endpmatrix
\quad \Bigg| x_{i}\in \bold{G}_{m}\right\}
$$
We denote a typical element of $\TT$ by $t(\{x_i\}).$
The simple roots are given by
$\{e_{i}- e_{i+1}\}^{n-1}_{i=1} \cup \{2e_{n}\},$ where
$e_{i}(t(\{x_i\}))= x_{i}, \ 2e_{i}(t)= x^{2}_{i}.$
This root system is of type $C_n.$
For the coroots, we note that $\check e_{i(x)} = (2e_{i})\spcheck(x)=
t(\{1,1,\dots,x,1,\dots,1\}),$
where the $x$ appears in the $i$\snug-th position.
Similarly,
$$(e_i-e_j)\spcheck(x)=t(\{1,\dots,x,\dots,x^{-1},1,\dots,1\}),$$
where $x$ appears in the $i$\snug-th position, and $x^{-1}$ is in the
$j$\snug-th position. So, $\check\Psi$ is of type $B_{n}.$
Therefore,
$$
^L\! G^0=SO(2n+1,\Bbb C)=\{g\in GL_{2n+1}(\Bbb C)|^t\! g J'g=J'\}
$$
$$
J'= \pmatrix 0&0&I_n\\0&1&0\\I_n&0&0\endpmatrix
$$
Since $\G$ is a split group,
$^L\! G=^L\! G^0\times \Gamma_{L/F}$

(3) $\G = SO_{2n+1}, \quad {}^L\! G= Sp_{2n}(\Bbb C)\times \Gamma_{L/F}$

(4) $\G = SO_{2n}, \quad {}^L\! G= Spin(2n)\times \Gamma_{L/F}$

(5) $\G = U_{n,n}.$  Let $E/F$
be a quadratic extension,  with $\s:x\mapsto \overline x$ the
Galois automorphism.
Let $E=F(\b),$ with $\bar\b=-\b.$  Set
$\quad J=\pmatrix 0&\b I\\
-\b I&0\endpmatrix.$ Then
$\G=\{g\in GL_{2n}|{}^t \! \overline g Jg = J\}.$
Looking at the maximal torus of diagonal
elements, we see that $\G$ is quasi-split, but not split.
One can see that
$\bold{G}(E)=GL_{2n}(E)$ and thus,
${}^L \! G^0=GL_{2n} (\bc).$

Let $x\in {}^L \! G^0.$  Then
$\sigma(x)=\Phi_n\ ^t x^{-1} \Phi_n^{-1},$
where
$$\ \Phi_n=\pmatrix &&&&&\cdot\\
&&&&\cdot\\
&&&\cdot\\
&&-1\\
&1\\
-1\endpmatrix$$
(see\Lspace \Lcitemark 17\Rcitemark \Rspace{}).  Then this gives the action of
$\Gamma_{E/F}$
on $^L\!G^0,$ and thus defines $^L\!G.$
\endexample

\subheading{(b) The  first conjecture}  The heart of the Langlands program is
the philosophy that harmonic analysis, number theory, and geometry should  be
connected within the theory of automorphic forms.  The formalization of this
philosophy lies in the very deep conjectures of Langlands.  We give
a rough sketch of one of these conjectures below.

Let $\G$ be as in Section 1, and set $G=\G(F).$  We denote by $P$
the projection of $^L\!G$ onto $\Gamma_{L/F}.$
 Let $\ph:\Gamma_{L/F}\lra \ ^L\!G$ be a homomorphism.  We say that $\ph$ is
{\bf admissible} if $P\circ\ph$ is the identity map.  We say that two
admissible homomorphisms are equivalent, if they differ by an inner
automorphism of $^L\!G.$  Notice that if $\G=GL_1=\G_m,$ then an admissible
homomorphism is just a character of $\Gamma_{L/F}.$

\proclaim{Conjecture 1.1 (Langlands)}
The equivalence classes of irreducible admissible representations of $G$
should be  parameterized by the equivalence classes of admissible
homomorphisms $\varphi\colon \Gamma_{L/F} @>>> {}^L\! G.$
(In fact we should replace $\Gamma_{L/F}$ by the Weil--Deligne group
$W_F'$ {\rm\Lspace \Lcitemark 7\Citecomma
80\Rcitemark \Rspace{}}. Of course, in this context, admissible is
a somewhat more technical concept.)
\endproclaim

\remark{Remarks}
 When $\G=\G_m=GL_1,$ we see that,
for a fixed $L,$ the admissible homomorphisms $\ph:\Gamma_{L/F}\lra\ ^L\!G$
parameterize only those characters which "factor" through $L.$
Thus, we must consider all fields $L$ over which $\G$ splits.  This is why
one must replace $\Gamma_{L/F}$ with a larger group.
However,  even with the Weil-Deligne group,
Conjecture 1 is too much to ask for.
For instance, for $\G=SL_2,$ Labesse and Langlands,
\Lcitemark 50\Rcitemark \Rspace{}, showed that sometimes
inequivalent representations must be parameterized
by the same admissible homomorphism.
So, one can only expect to partition the
tempered representations of $G$  into finite subsets, called $L$-{\bf packets},
such that these $L$-packets are parameterized by admissible
homomorphisms, modulo conjugacy in ${}^L\! G.$  The representations in a
given $L$-packet are said to be $L$-{\bf indistinguishable}.
\endremark

\example{Example} If $\pi$ is a tempered representation of $GL_{n}(F),$
then $\{\pi\}$ is an $L$-packet\Lspace \Lcitemark 7\Rcitemark \Rspace{}.
Suppose $\pi|_{SL_{n}(F)}=\pi_{1}\oplus \ldots \oplus \pi_{k}.$
Then\Lspace \Lcitemark 78\Rcitemark \Rspace{} the $\pi_i$ are distinct and
Gelbart and Knapp,\Lspace \Lcitemark 21\Rcitemark \Rspace{},
showed that $\{\pi_1,\ldots,\pi_k\}$ is an
$L$-packet for $SL_{n}(F).$
Note $^L\! (SL_n(F))^0=PGL_n(\Bbb C).$
Suppose $\varphi:\Gamma_{L/F}\lra GL_n(\bc)\times\Gamma_{L/F}$
is a parameter for $\pi.$
Composing with the projection $\eta:GL_n(\bc)\lra PGL_n(\bc),$
should give a parameter for an $L$\snug-packet of $SL_n(F).$
What Gelbart and Knapp showed was that, assuming
the Langlands correspondence is understood for $GL_n,$
then $\eta\circ\varphi$ must be the parameter for $\{\pi_1,\dots,\pi_k\}.$
\endexample

\subheading{(c) Unramified representations} Since $\bold{G}$
splits over an unramified extension, $L$ of $F,$ we can take the $\Cal
O_{F}$-points of $\G.$
Let $K= \bold{G}(\Cal O_{F}).$  Then, $K$
is a ``good'' maximal compact subgroup of $G$\Lspace \Lcitemark 11\Rcitemark
\Rspace{}.  A
representation $(\pi, V)$ of $G$ is {\bf unramified}, or
{\bf class 1}, if there is a $v\in V$ with $\pi(k)v=v$ for all
$k\in K.$

\proclaim{Lemma 1.2} (See {\rm\Lspace \Lcitemark 11\Rcitemark \Rspace{}}.)
If $(\pi,V)$ is an admissible irreducible representation of $G=\G(F),$
then
$\dim_\bc(V^K)\leq 1.$
\endproclaim

Suppose that $\pi$ is class 1.
Note that $\left(\pi\left(\Cal H\left(G//K\right)\right), \ V^K)\right)$
is a character, $\chi,$ and $f\mapsto
\chi(f)$ determines a semi-simple conjugacy class $\{A\}$ in ${}^L\! T^0,$ via
the
Satake isomorphism\Lspace \Lcitemark 14\Rcitemark \Rspace{}, unique up to the
action of the Weyl group $W(^L \! G^0, ^L \! T^0).$

\example{Example}Suppose $\G=GL_{n}$
and $\pi = \Ind_{B}^G(\omega_{1},\ldots \omega_{n}),$ with $\omega_{1},\ldots
\omega_{n}$ unramified characters,
then
$$
A=\pmatrix \omega_{1}(\varpi)&&\\
&\ddots&\\
&&\omega_{n}(\varpi)\endpmatrix,
$$
and $A$ determines $\pi$ (up to permutation of $(\omega_{1},\ldots,
\omega_{n})).$
\endexample

Let $\tau$ be the Frobenius class in $\Gamma_{L/F}.$ Suppose that $\pi$ is
unramified.
Suppose r is a
finite dimensional representation of ${}^L\! G,$ i.e., $r: {}^L\! G\to
GL_{n}(\Bbb C)$ is a homomorphism, with $r|_{{}^L\! G^{0}}$ a complex
analytic representation of ${}^L\! G^{0}.$  Let $\widetilde r$ be its
contragredient.  For $s\in\Bbb C,$ let
$$
L(s,\pi,r)= \det(I- r(A\rtimes\tau)q^{-s})^{-1}.
$$

In trying to understand the reason for assigning
this value to $\pi,\ r,$ and $s,$
one should keep in mind the concept of an Artin $L$-function
\Lcitemark 19\Citecomma
32\Rcitemark \Rspace{}.
If we have a representation $r\colon \Gamma_{L/F}\to
GL_{n}(\bc),$ then
the attached local Artin $L$-function is $L(s,r)=\det(I-r(\tau)q^{-s})^{-1}$
and we hope this determines the splitting of the prime ideal in $L.$

\subheading{(d) Automorphic Representations on $G$}
Let $F$ be a number field,
and let $\Bbb A_{F}$ be the adeles of $F.$
Let $\pi=\bigotimes\limits_v\pi_v$ be a cuspidal automorphic
representation
of $\G(\BA_F).$  We refer
to\Lspace \Lcitemark 8\Citecomma
54\Rcitemark \Rspace{} for the definition
of the space
$L^2_0(\G(F)\backslash \G(\BA_F),\omega)$
of cuspidal representations whose central character is $\omega.$

Note that
$\G_v=\bold{G}\underset F\to\times F_{v}$ gives the group $\bold{G}$
over $F_{v},$ and we have ${}^L\! G$ ($L$-group of $\G$) and
${}^L\! G_{v},$ the $L$--group of $\G_v.$
For almost all $v, \ \bold{G}\underset F\to\times F_{v}$ is unramified, and
$\pi_{v}$ is $K_{v}$ unramified.  Suppose $r$ is a representation of
${}^L\! G,$ and let $r_{v}$ be defined by $r_v: {}^L\! G_{v}\to
{}^L\! G\overset r\to\longrightarrow GL_{n}(\Bbb C).$  Then we have a local
$L$-function, $L(s,\pi_{v}, r_{v}),$ whenever $\pi_{v}$ and
$\bold{G}_{v}$ are unramified.

\subheading{(e) The Main Conjecture} We now give Langlands conjecture
on the existence of global $L$\snug-functions
\Lcitemark 51\Citecomma
52\Rcitemark \Rspace{}.
Let $\psi$ be a character of $\Bbb A_F$ which is trivial
on $F,$ and suppose $\psi=\prod\limits_v\psi_v.$
Let $S$ be a finite set of places so that if $v\notin S, \ \pi_{v}$ and
$G \underset F\to\times F_{v}$ are unramified.  Let
$$
L_{S}(s,\pi,r)= \underset{v\notin S}\to \prod L(s,\pi_{v}, r_{v}).
$$

\proclaim{Theorem 1.3} ({\rm\Lspace \Lcitemark 7\Rcitemark \Rspace{}}).
$L_{S}(s,\pi,r)$ converges for
$Re\ s > > 0.$
\endproclaim

\proclaim{Conjecture 1.4 (Langlands)}
For $v\in S$ it is possible to define a local $L$-function
$L(s,\pi_{v}, r_{v}),$ so that $L(s,\pi_{v}, r_{v})=(P_{v}
(q^{-s}_{v}))^{-1},$ with $P_{v}(t)$ a polynomial whose constant
term is $1,$ and a local root number
$\varepsilon(s, \pi_{v}, r_{v}, \psi_{v}),$ (a monomial in
$q^{-s}_{v})$ so that
$$
L(s,\pi, r) = \underset v\to\prod L(s,\pi_v,r_v)
$$
has a meromorphic continuation to $\Bbb C,$ with finitely many poles, and
$$
L(s,\pi,r)=\varepsilon(s,\pi,r) L(1-s, \pi,\widetilde r),$$
with
$$\varepsilon(s,\pi,r)=\underset v\to\prod \varepsilon(s,\pi_{v}, r_{v},
\psi_{v})$$
Moreover, if $v\in S,$ then $\varepsilon(s,\pi_{v},r_v,\psi_{v})=1,$ and
$L(s,\pi_v,r_v)$ is a s in subsection {\bf (d)}.
\endproclaim

\remark{Remarks}
\roster
\item Note that the local root numbers $\e(s,\pi_v,r_v,\psi_v)$ depend on
the choice of the character $\psi,$ but that the global root number
$\e(\pi,s,r)$ does not.

\item  Suppose that $\varphi_v:\Gamma_{\bar F_v/F_v}\lra \phantom{}^LG_v$
is the admissible homomorphism corresponding to $\pi_v.$  Then
$L(s,\pi_{v}, r_{v})$ should be $L(s,r_{v}\circ\varphi_{v}),$
where this last object is the Artin $L$--function.

\item  One must be careful not to read too much into this conjecture.
In particular, it does not formally describe the nature of a global
Langlands parameter.  Such a parameterization requires an object
(not yet known) much larger than the global Weil-Deligne group.

\item
Because of the relationship between automorphic forms
and automorphic representations,
$L$-functions, which are arithmetic in nature, must play a fundamental role
in the harmonic analysis of reductive groups over both local and
global fields.

\endroster
\endremark

While this exposition of the Langlands program is far from complete, we hope
that it will give the reader enough background to see how the theory of the
next section fits into this philosophy.

\subheading{\S 2  Intertwining operators and reducibility of induced
representations}

\subheading{(a) Preliminaries}
If $X$ is a totally
disconnected space, and $Y$ is a complex vector space, then we let
$C^\infty(X,Y)$ be the space of functions $f:X\lra Y$
which are locally constant.  We let
$C_c^\infty(X,Y)$ be the subspace of $f\in C^\infty(X,Y)$
which are compactly supported.
For a totally disconnected group $G,$
we let $\Cal E_c(G)$ be the collection of equivalence classes
of irreducible admissible representations of G.
We denote by $\Cal E(G)$ the  (pre)-unitary classes in $\Cal E_c(G).$
The supercuspidal classes are denoted by $^\circ\snug{\Cal E}_c(G),$
and $^\circ\Cal E(G)$ denotes  $\Cal E(G)\cap\ ^\circ\Cal E_c(G).$
We let $\Cal E_2(G)$ denote the discrete series, and $\Cal E_t(G)$ the
tempered classes in $\Cal E_c(G).$  (See\Lspace \Lcitemark 11\Rcitemark
\Rspace{}
for details on these definitions.)

Let $F$ be a local field of characteristic zero.  Suppose $\bold{G}$ is a
connected reductive, quasi--split algebraic group, defined over $F.$  Let
$G=\bold{G}(F).$
Suppose $\bold{P}$ is a parabolic subgroup of $\bold{G},$  and let
$\bold{P}=\bold{M}\bold{N}$ be the Levi decomposition
of $\PPP.$  We call $\M$ the Levi component of $\PPP,$
and $\N$ the unipotent radical of $\PPP.$
Let $P=\PPP(F)=MN=\M(F)\N(F).$
Suppose that $(\s,V)$ is an admissible complex representation of $M.$
We denote the contragredient representation by
$(\wt \s,\wt V).$
We let
$\d_P$ be the modular function of $P.$  With this data, we set
$$V(\s)=\{f\in C^\infty(G,V)|f(mng)=\s(m)\d_P(m)^{1/2}f(g),\ \forall m\in M,
n\in N, g\in G\}.$$
Then $G$ acts on $V(\s)$ by right translation, and this action
is called the representation of $G$ {\bf unitarily induced} from
$\s.$  We denote the induced representation
by $\Ind_P^G(\s).$   The factor $\d_P^{1/2}$ is
there to ensure that $\Ind_P^G(\s)$ is unitary  if
$\s$ is. (Hence the term unitary induction.)

The classification theorems of Jacquet,\Lspace \Lcitemark 31\Rcitemark
\Rspace{},
and Langlands,\Lspace \Lcitemark 9\Rcitemark \Rspace{}, indicate the importance
of studying these induced representations.  When $\s$
is an irreducible discrete series representation,
then the components of $\Ind_P^G(\s)$ are tempered.
Furthermore,  every irreducible tempered representation of $G$
is a component of $\Ind_P^G(\s)$ for some $P,$ and some
discrete series $\s$  of $M.$   In this section we are interested in
two aspects of induced representations.  The first is the classification
of the tempered spectrum of $G,$ i.e., determining the structure
of $\Ind_P^G(\s)$ for every choice of $\PPP=\M\N,$
and every irreducible discrete series representation $\s.$
That is, we wish to determine for which $\s$ the representation
$\Ind_P^G(\s)$ is reducible.
Furthermore, we want to know how many components
$\Ind_P^G(\s)$ has, what are the multiplicities with which these components
appear, and how are the characters of these components related.
The second point of interest is determining the arithmetic properties
of $\Ind_P^G(\s),$ and its components.  In particular
one wants to know how the $L$\snug-functions for $\s$
and those for the components of $\Ind_P^G(\s)$ are related. The fact that
the answers to these two questions are related is quite deep,\Lspace \Lcitemark
67\Rcitemark \Rspace{},
and is in fact what has allowed significant progress to be made
recently.

Suppose that $\ph:W_F\lra\ ^LM$ is a conjectural parameter for the discrete
series $L$\snug-packet $\{\s\}$ which contains $\s.$
Then
$W_F\overset\ph\to\lra\ ^LM\overset{i}\to\hookrightarrow\ ^LG,$
should define an $L$\snug-packet for $G.$  If
$\Pi_G(\s)$ is the collection of components of
$\Ind_P^G(\s),$ then $\Pi_\ph(G)=\bigcup\limits_{\tau\in\{\s\}}\Pi_G(\tau),$
should be this $L$\snug-packet.  Given the conjectural properties
of the Langlands $L$\snug-functions given in Section 1,
it makes heuristic sense to believe that the $L$\snug-functions
for elements of $\Pi_\ph(G)$ are ``induced" from
those of $\{\s\}.$  That is, for $\pi\in\Pi_\ph(G),$
and a complex representation $\rho$ of $^LG,$
we expect that $L(\nu,\pi,\rho)=L(\nu,\rho\circ i\circ\ph),$
where the last object is the Artin $L$\snug-function.

Let $\A$ be the split component of $\M,$ i.e., maximal
torus in the center of $\M.$
Then $\bold{M}=Z_\bold{G}(\bold{A}).$
Let $W=N_{\bold{G}}(\bold{A})/\bold{M}.$
Then $W$ is called the {\bf Weyl group} of $\G$ with respect to $\A,$
and we may denote this group by $W(\G,\A),$ if there is any
ambiguity about $\G$ and $\A.$

Let
$\wt w\in W.$
Choose $w\in N_{G}(A)$ representing $\wt w.$
Let $w\sigma(m)=\sigma(w^{-1}m w).$
The class of $w\sigma$
is independent of the choice of representative $w$ for $\wt w.$
We write $W(\s)$ for the subgroup of $\wt w\in W$ which fix the class of $\s.$
Let $C(\s)$ be the commuting algebra of $\Ind_P^G(\s).$

\proclaim{ Theorem 2.1 (Bruhat\Lspace \Lcitemark 10\Rcitemark \Rspace{})}
For  $\s\in\Cal E_2(M),$ we have $\dim_\bc(C(\s))\leq |W(\s)|.$
\qed
\endproclaim

One wishes to use the group $W(\s)$ to decompose $\Ind_P^G(\s).$
This led to the development of the theory of the standard
intertwining operators.  Namely, one can attach to each element
$w$ of $W(\s)$ a self intertwining operator $\Cal A(w,\s)$
for $\Ind_P^G(\s),$ and a complete
understanding of these operators determines the algebra $C(\s).$
We will outline the theory of these operators.

Set $\frak a=\text{Hom}(X(\M)_{F}, \Bbb R),$ where
$X(\M)_{F}$ is the set of $F$--rational characters of $\M.$
Let $\frak a^*=X(\M)_F \otimes_{\Bbb Z} {\Bbb R}.$
Let $\frak a^{*}_{\Bbb C}=\frak a^{*}\bigotimes\limits_{\Bbb R}\Bbb C.$
Then $\frak a^*_{\Bbb C}$ is a complex manifold.
There is a homomorphism, $H_P\colon X(\bold{M})_{F}\to \frak a,$
so that, $|\chi(m)|_{F}=q^{\langle\chi, H_{P}(m)\rangle},$ for all
$\chi\in X(\bold{M})_{F}, \ m\in\M$\Lspace \Lcitemark 31\Rcitemark \Rspace{}.
Let $\nu \in \frak a^{*}_{\Bbb C},$
and $I(\nu,\sigma)=\text{Ind}^{G}_{P}(\sigma\otimes q^{\langle\nu,
H_{P}()\rangle}).$
We denote by
$V(\nu,\s)$ be the space of functions on which we realize
$I(\nu,\s).$

Let $K$ be a good maximal compact (say  $G(\Cal O_{F})$).
(See\Lspace \Lcitemark 12\Citecomma
56\Rcitemark \Rspace{} for a definition
of a good maximal compact.)
Then, for $k\in M\cap K,$ and all $\nu,$
we have $\langle\nu, H_{P}(k)\rangle=1.$
Thus, we suppose  that $f\in C^\infty(K,V)$ is such that
$$
f(m_{K}n_{K},k^{\prime})=\sigma(m_{K})\delta^{1/2}_{P}(m_{K})f(k^{\prime}),\text{ for all } m_Kn_K\in P\cap K,\ k'\in K.
$$
By the Iwasawa decomposition,\Lspace \Lcitemark 11\Citecomma
56\Rcitemark \Rspace{}, $G = PK.$  So, for
each $\nu,$ we set $f_{\nu}(mnk) = \sigma(m)q^{\langle\nu,
H_{P}(m)\rangle}\delta_P^{1/2}(m)f(k)\in V(\nu,\sigma).$
Note $f=f_{\nu}|_{K}.$
Thus, we have a natural isomorphism between $V(\nu,\sigma)$
and $V(\nu^{\prime},\sigma)$ given by $f_{\nu'}\leftrightarrow f_{\nu}.$

Fix $\nu\in\frak a^{*}_{\Bbb C},$ and $\sigma\in \Cal E_{2}(M).$
Let $\wt w\in W,$ and choose $w$ representing $\wt w.$
Fix $f\in I(\nu,\sigma).$
Define
$$
A(\nu,\sigma,w)f(g) = \int\limits_{N_{\wt w}}f(w^{-1}ng)\, dn,
$$
where $\N_{\wt w}= \U\cap w^{-1} \N^- w,$
and $\N^-$ is opposite to $\N.$
Let $m_{1} \in M, \ n_{1}\in N.$  Then,
$$
\gathered
A(\nu,\sigma,w)f(m_1 n_1 g)= \int\limits_{N_{\wt w}}f( w^{-1}nm_1 n_1 g)\, dn\\
=\int\limits_{N_{\wt w}}f( w^{-1}m_1  w w^{-1}m_{1}^{-1}n
m_1 n_1 g)\, dn\\
=\sigma( w^{-1}m_1  w)\delta^{1/2}_{P}
({w}^{-1}m_1  w)
q^{\langle\nu, H_{P}({w}^{-1}m_1 {w})\rangle}\int\limits_{N_{\wt w}}
f({w}^{-1}m_1^{-1}nm_1n_1g)\, dn.
\endgathered$$

\noindent
Now we note that, on $N_{\wt w},$ the measures $d(m_1^{-1}nm_1),$
and $dn$ are related by
$$\aligned
d(m_1^{-1}nm_1)&=\left({\d_P(m_1)}/{w\d_P(m_1)}\right)^{1/2}\\
&=\left({\d_P(m_1)}/{\d_P(w^{-1}m_1w)}\right)^{1/2},
\endaligned$$
(cf\Lspace \Lcitemark 62\Rcitemark \Rspace{}).
Therefore,
making the substitution
$n'=m_1^{-1}nm_1n,$ we have
$$
\aligned
A(\nu,\sigma,w)f(m_1 n_1 g)=& w\sigma(m_1)\delta^{1/2}_P(m_1)q^{\langle w\nu,
H_{P}(m_1)\rangle}\int\limits_{N_{\wt w}}
f({w}^{-1}n'g)\, dn'\\
=& w\sigma(m_1)\delta^{1/2}_P(m_1)q^{\langle w\nu, H_{P}(m_1)\rangle}
A(\nu,\sigma,w)f(g),\\
\endaligned
$$
which implies
$A(\nu,\s,w)f\in V(w\nu, w\sigma).$

Note that we have only defined $A(\nu,\s,w)$ formally, in
that we have said nothing about the convergence of $A(\nu,\s,w)f(g).$
However, the above argument shows that, if $A(\nu,\s,f)(g)$
converges for all $f\in V(\nu,\s),$ and $g\in G,$ then
$A(\nu,\s,w)$ defines an intertwining operator between $I(\nu,\s)$
and $I(w\nu,w\s).$

\proclaim{Theorem 2.2 (Harish--Chandra)} Let $\sigma\in \CalE_2 (M).$
Then, for $\text{Re }\nu > 0,$ the operator $A(\nu,\sigma,w)$ converges
absolutely, with
a meromorphic continuation to all of $\frak a^*_{\bc}.$
So,  fixing $f\mapsto f_\nu,$ as above, then
$\nu \mapsto \langle \widetilde{v}, A(\nu,\sigma,w) f_\nu (g)\rangle$ is
meromorphic for all $f,\widetilde{v},$ and $g.$\qed
\endproclaim

Suppose $w\sigma\simeq \sigma.$
Then $A(0,\sigma,w)$ gives an intertwining operator between
$\Ind_P^G(\s)$ and $\Ind_P^G(w\s),$  which are isomorphic
representations.
Of course, $A(\nu,\sigma,w)$ may have a pole at $\nu=0,$
and this question of analyticity at $\nu=0$ is in fact
crucial to determining the structure of $\Ind_P^G(\s).$

\proclaim{Theorem 2.3 (Harish--Chandra)} There is a meromorphic, complex valued
function $\nu \mapsto \mu (\nu,\sigma,\wt w)$ so that
$$
A(\nu,\sigma,w) A(w\nu,w\sigma,w^{-1})=\gamma_{\wt w}
(G/P)^2 \mu (\nu,\sigma,\wt w)^{-1},
$$
where
$${\gamma_{\wt w} (G/P)=\int\limits_{N_{\wt w}} q^{\langle \nu, H_P (n)\rangle}
\, dn}.$$
Furthermore, $\mu(\nu,\sigma,w)$ is holomorphic and non--negative on
$i\frak a^*.$\qed
\endproclaim

Let $\Phi(\PPP,\A)$ be the reduced roots of $A$ in $P.$
The {\bf length} of $\wt w\in W,$ is given by $\ell(\wt w)=|\{\alpha\in \Phi
(P,A)|\wt w\alpha < 0\}|.$
There is a longest element $\wt w_0\in W$\Lspace \Lcitemark 15\Rcitemark
\Rspace{}.
We write $\mu(\nu,\sigma)$ for $\mu(\nu,\sigma,\wt w_0),$ and $\mu(\sigma)$
for $\mu(0,\sigma).$  We call
$\mu(\nu,\sigma)$ the Plancherel measure of $(\nu,\sigma).$

\proclaim{Theorem 2.4 (Harish--Chandra\Lspace \Lcitemark 72\Rcitemark
\Rspace{}})
Suppose $P$ is maximal and proper, and there is some $\wt w\not= 1$ in $W,$
with
$\wt w\sigma\simeq\sigma.$
Then $\Ind_P^G(\sigma)$ is reducible if and only if $\mu(\sigma)\not= 0.$\qed
\endproclaim

\def\lp{^L\! P}

\def\lm{^L\! M}
\def\ln{^L\! N}

\subheading{(b) Arithmetic Considerations}

We now examine the arithmetic properties of the function $\mu.$  It turns out
that the Plancherel measure is directly related to the theory of
$L$\snug-functions.  Suppose, for the moment that $\PPP$ is a maximal proper
parabolic subgroup of $\G.$
Then $N_{\wt w}=N,$ and $\frak a^*_\bc/\frak z\simeq\bc.$
(Here $\frak z$ is the Lie algebra of the split component
of $\G.$)
Let $\rho$ be a representation of $\lm,$ and set
$$
\gamma (s,\sigma,\rho,\psi_F)=\varepsilon (s,\sigma,\rho,\psi_F)
L(1-s,\widetilde\sigma,\rho)/L(s,\sigma,\rho),
$$
where $L(s,\s,\rho)$ is the conjectural Langlands
$L$\snug-function attached to $\s$ and $\rho.$
By\Lspace \Lcitemark 7\Rcitemark \Rspace{} there is a parabolic subgroup $^L\!
P$ of $^L\! G,$ with
$\lp= {}^L\! M \ln,$ for some unipotent subgroup $\ln.$
 Let $^L\! \frak{n}$ be the Lie algebra of $\ln.$
Denote by $r$ the adjoint representation of $\lm$ on $^L\! \frak{n}.$

\proclaim{Conjecture 2.5 (Langlands\Lspace \Lcitemark 53\Rcitemark \Rspace{})}

If $\s\in\Cal E_2(M),$ then
$\mu(\sigma) \gamma_{\wt w_0}(G/P)^2=
\gamma(0,\sigma,r,\overline{\psi}_F)\gamma(0,\widetilde\sigma,r,\psi_F).$\qed
\endproclaim

To proceed, we need the notion of a generic representation.
Suppose that $\G$ is a connected reductive
quasi-split algebraic group, defined over $F.$
Let $\Phi$ be the roots of $\bold T$ in $\G,$
and $\DE$ the set of simple roots given by our choice
of a Borel subgroup.
For $\a\in\DE,$ let $\U_\a$ be the subgroup of $\U$
whose Lie algebra is $\g_\a\oplus\g_{2\a}$\Lspace \Lcitemark 6\Citecomma
35\Citecomma
57\Rcitemark \Rspace{}.
The subgroup
$$\align
U'&=\prod_{\a\in\Phi^+\setminus\DE}U_\a\\
\intertext{is normal in $U.$   Furthermore,}
U/U'&\simeq\prod_{\a\in\DE}U_{\a}/U_{2\a}.
\endalign$$
Let $\a\in\DE$ and let $\psi_\a$
be a character of $U_{\a}/U_{2\a}.$
The character $\psi$ of $U,$ trivial on $U',$
and given by
$$\psi=\prod_{\a\in\DE}\psi_\a,$$
is called {\bf non-degenerate} if $\p_\a$ is non-trivial
for each $\a.$  An admissible representation
$(\pi,V)$ of $G$ is called {\bf non-degenerate,} or
{\bf generic} if there is a  non-degenerate
character $\p$ of U, and a linear functional $\l$ on
$V,$ such that
$$\l(\pi(u)v)={\p(u)}\l(v),\ \text{ for all }\  u\in U, v\in V.$$
Such a functional is called a {\bf Whittaker functional.}
Let $V_\p^*$ be the complex vector space of Whittaker functionals
on $V$ (with respect to a fixed $\p$).

\example{Example}  Let $\G=GL_n,$ so $\G(F)=GL_n(F).$
Suppose that $\B=\TT\U$ is the Borel subgroup of upper triangular
matrices.  For $t\in F,$
we denote by $E_{ij}(t)$ the $n\times n$ matrix whose
$ij$\snug-th entry is $t,$ and all other entries are zero.
Let $\a=e_i-e_{i+1}.$  Then
$U_\a=\{I+E_{i(i+1)}(t)|t\in F\}.$ Recall that $\psi_F$ is a fixed
non-trivial character of $F^+.$  Let $\psi_\a(I+E_{i(i+1)}(t))=\psi_F(t).$
Setting $\psi=\prod\limits_{\a\in\DE}\p_\a,$
we have
$$\p\left(\pmatrix
1&x_{12}\\
&1&x_{23}&&*\\
&&\ddots\\
&&&1&x_{(n-1)n}\\
&&&&1\endpmatrix\right)=\p_F(x_{12}+x_{13}+\dots+x_{(n-1)n}).$$
Then $\p$ is non-degenerate, and in fact any
non-degenerate character of $U$ is conjugate to $\p.$

It is a result of Gelfand and Kazhdan that every irreducible
discrete series representation of $GL_n(F)$ is generic\Lspace \Lcitemark
23\Rcitemark \Rspace{}.
Jacquet,\Lspace \Lcitemark 37\Rcitemark \Rspace{}, showed that every
irreducible
tempered representation of $GL_n(F)$ is generic.
For other classical groups, there are examples of discrete
series representations which fail to be generic\Lspace \Lcitemark 34\Rcitemark
\Rspace{}.
\endexample
\proclaim{Theorem 2.6 (Shalika\Lspace \Lcitemark 70\Rcitemark \Rspace{})}
Let
$(\pi,V)$ be  an irreducible admissible representation
of $G,$ then $\dim_\bc V_\psi^*\leq 1.$\qed
\endproclaim

Rodier showed that the dimension of $V_\psi^*$ is preserved
under parabolic induction, which we state below.

\proclaim{Theorem 2.7 (Rodier\Lspace \Lcitemark 59\Rcitemark \Rspace{})}
Let $\G$ be quasi-split, and suppose that $\PPP=\M\N$ is a
parabolic subgroup of $\G.$  If $(\s,V)$ is an
irreducible admissible representation of $M,$
and $(\pi,W)\simeq \Ind_P^G(\s),$  then
$\dim_\bc W_\psi^*=\dim_\bc V_\psi^*.$\qed
\endproclaim

Suppose $(\s,V)$ is an irreducible generic discrete series representation
of $M.$  Let $\l\in V_\psi^*$ be non-zero.  For each $v\in V,$
let $W_v(m)=\l(\s(m)v).$  Then $W_v$ is called the
{\bf Whittaker function} attached to $\l$ and $v.$
Note that, for $u\in U\cap M,$ we have
$W_v(ug)=\l(\s(ug)v)=\p(u)W_v(g).$  Thus, $W_\s=\{W_v|v\in V\}$
is a subspace of $\Ind_{U\cap M}^M(\psi).$  For
$m_1\in M,$ we see that
$$m_1\cdot W_v(m)=W_v(mm_1)=\l(\s(mm_1)v)=W_{\s(m_1)v}(m).$$
So the restriction of $\Ind_{U\cap M}^M(\psi)$ to $W_\s$ is
isomorphic to $\s,$ and we call this the {\bf Whittaker model}
for $\s.$

Let $\nu\in\frak a_\bc^*.$  Consider $V(\nu,\s)$
as a space of functions from $G$ to $W_\s.$
For each $f\in V(\nu,\s),$ and each $g\in G,$
we denote by $(f(g),m)$ the value of $f(g)$ at $m.$
There is a natural way to define a Whittaker functional
$\l(\nu,\s,\psi)$ on $V(\nu,\s).$
Namely, let
$$\l(\nu,\s,\p)(f)=\int\limits_{N'}f(w_0^{-1} n',e)\p( n')\ dn',$$
where $\N'=w_0\N^- w_0^{-1}.$
Casselman and Shalika,\Lspace \Lcitemark 16\Rcitemark \Rspace{},
showed that $\l(\nu,\s,\p),$ which converges
absolutely in a right half plane,
has an analytic continuation to all of $\frak a_{\bC}^*.$
Moreover,
$\l(\nu,\s,\p)$ defines an element of $V(\nu,\s)_\p^*.$
Note that if $f\in V(\nu,\s),$ then, for each $u\in N,$
$$\gather
\l(w\nu,w\s,\psi)A(\nu,\s,w)(I(\nu,\s)(u)f)=\l(w\nu,w\s,\p)I(w\nu,w\s)(u)A(\nu,\s,w)f\\
=\p(u)\l(w\nu,w\s,\p)A(\nu,\s,w)f,
\endgather$$
so $\l(w\nu,w\s,\psi)A(\nu,\s,w)\in V(\nu,\s)_\p^*.$
Rodier's Theorem gives rise to the following result.

\proclaim{Proposition 2.8\ (\Lcitemark 64\Rcitemark )}
There is a complex number $C_\psi(\nu,\s,w),$ such that
$$\l(\nu,\s,\psi)=C_\psi(\nu,\s,w)\l(w\nu,w\s,\psi) A(\nu,\s,w).$$
Furthermore, $\nu\mapsto C_\psi(\nu,\s,w)$ is
meromorphic on $\frak a_\bc^*.$\qed
\endproclaim

We call $C_\psi(\nu,\s,w)$ the {\bf local coefficient} attached to
$\p,\nu,\s,$
and $w.$

\proclaim{Proposition 2.9\ (\Lcitemark 64\Rcitemark )}

\roster
\item For all $\nu\in\frak a_\bc^*,$ we have the identity
$$C_\psi(w\nu,w\s,w^{-1})=\overline{C_\psi(-\bar\nu,\s,w)}.$$
\item If $\nu = -\bar\nu,$ and
$\s$ is unitary, then
$$|C_\psi(\nu,\s,w)|^2=\gamma^{-2}(G/P)\mu(\nu,\s).\qed$$
\endroster
\endproclaim

By studying the relationship of local coefficients to Plancherel
measures, Shahidi was able to prove Langlands's conjecture on
Plancherel measures for generic representations.
He has also used this method to derive
estimates toward the Ramanujan conjecture\Lspace \Lcitemark 66\Rcitemark
\Rspace{}.

\proclaim{Theorem 2.10 (Shahidi\Lspace \Lcitemark 67\LIcitemark{}, Theorem
3.5\RIcitemark \Rcitemark \Rspace{})}
Let  $\sigma$  be irreducible admissible and generic.
Then, for each $r$ and $\psi_F,$
there exists a unique arithmetic function, $\gamma(s,\s,r,\psi_F),$
(cf Theorem 3.5 of {\rm\Lspace \Lcitemark 67\Rcitemark \Rspace{}} for its
uniquely defining properties) such that
$$\mu(\sigma) \gamma_{w_0}(G/P)^2=
\gamma(0,\sigma,r,\overline{\psi}_F)\gamma(0,\widetilde\sigma,r,\psi_F).$$
Moreover, if we accept two conjectures in harmonic analysis, we
can remove the assumption that $\s$ is generic.\qed
\endproclaim

\subheading{(c) $R$-groups}
We now move toward our goal of describing the intertwining algebra
$C(\s)$ of $\Ind_P^G(\s).$  The idea is to use Plancherel
measures, and the operators $A(\nu,\s,w),$ to construct
normalized intertwining operators which have no poles on
the unitary axis $i\afa^*.$
{}From these we construct a collection of self-intertwining operators
for $\Ind_P^G(\s),$ and determine how to construct a basis
of $C(\s)$ from among these.
We describe this below.

\proclaim{Theorem 2.11}
There is a meromorphic normalizing factor, $r(\nu,\sigma,w),$ so
that
$$\calA (\nu,\sigma,w)=r(\nu,\sigma,w)A(\nu,\sigma,w)$$
is holomorphic on $i\frak a^*.$\qed
\endproclaim

\proclaim{Theorem 2.12 (Shahidi)\Lspace \Lcitemark 67\Rcitemark \Rspace{}} Let
$\sigma$ be
irreducible admissible and generic.  Then $r(\nu,\s,w)$ can be
chosen so that
$$
r(0,\sigma,w)=\varepsilon(0,\sigma,\widetilde{r},\psi_F)
L(1,\sigma,\widetilde{r})/ L(0,\sigma,\widetilde{r}).\qed
$$
\endproclaim

Now we see that reducibility of induced representations has an
arithmetic interpretation.
Conversely, if we can determine when $\Ind_P^G(\sigma)$ is reducible,
then we
can determine the poles of $L(s,\sigma,r).$

Let $A(\sigma,w)=A(0,\sigma,w),$ and $\calA(\sigma,w)=\calA(0,\sigma,w).$

\proclaim{Proposition 2.13} Let $w_1,w_2\in W$
\item{(1)}If $\ell(w_1w_2)=\ell(w_1)+\ell(w_2),$ then
$$
A(\sigma,w_1 w_2)=A(w_2\sigma,w_1)A(\sigma,w_2)
$$

\item{(2)}$\calA(\sigma,w_1 w_2)=\calA(w_2 \sigma,w_1)\calA(\sigma,w_2).\qed$
\endproclaim
We call (2) the cocycle relation for the normalized intertwining operators
$\Cal A(\s,w).$

We wish to determine the structure of $\Ind_P^G(\s)$
for arbitrary $P.$  So far we have Harish-Chandra's theorem
for maximal parabolic subgroups.  We want to know how to
utilize this theorem in the more general case.
Harish-Chandra proved a product formula for Plancherel
measures, which reduces the computation of Plancherel measure
to the case of maximal parabolic subgroups.

Let $\Phi(\PPP,\A)$ be the reduced roots of $\A$ in $\PPP.$
For $\beta\in \P(\PPP,\A),$ let
$\A_\beta=(\A\cap \ker \chi_\beta)^0$ and $\M_\beta=Z_{\G}(\A_\beta).$
(Here $\chi_\beta$ is the root character attached to $\beta.$)
Let $\N_\beta=\M_\beta\cap \N.$
Then ${}^*\PPP_\beta=\M \N_\beta$ is a maximal parabolic subgroup of
$\M_\beta.$
So, there is a Plancherel measure $\mu_\beta (\nu,\sigma)$ attached to
$\beta,\nu,$ and $\sigma.$  Note that $\mu_\b(\s)=0$ if and only if
$W(\M_\b,\A)\cap W(\s)\neq\{1\}$ and $\Ind_{^*\!{P_\b}}^{M_\b}(\s)$ is
irreducible.

\proclaim{Theorem 2.14 (Harish--Chandra\Lspace \Lcitemark 31\Rcitemark
\Rspace{}) Product formula for Plancherel Measures}
$$
\gamma^{-2} (G/P)
\mu(\nu,\sigma)=\prod_{\beta\in\Phi(\PPP,\A)}\gamma_\beta^{-2}
(M_\beta/^*P_\beta)\mu_\beta(\nu,\sigma).\qed
$$
\endproclaim

Notice that this gives us an inductive formula for
the $\gamma$\snug-factors, $\gamma(\nu,\s,w),$
which is one of their fundamental (conjectural) properties.
That such an inductive formula for $\gamma$\snug-factors
exists is part of Shahidi's result\Lspace \Lcitemark 67\Rcitemark \Rspace{}.

For $\b\in \P(\PPP,\A),$ let $\wt w_\b$ be the reflection in the reduced root
$\b.$  If $\wt w_\b\in W(\s),$ then $N_{\wt w_\b}=N_\b,$
and thus
$$A(\nu,\s,w_\b)f(g)=\int\limits_{N_\b}f(w_\b^{-1}ng)\, dn,\tag 2.1$$
which, for $g\in M_\b$ is the intertwining operator that
determines the reducibility of $\Ind_{^*P_\b}^{M_\b}(\s).$
Furthermore, the product formula for the intertwining operators
shows that every $A(\nu,\s,w)$ can be written as the composition
of operators of the form (2.1)\Lspace \Lcitemark 64\Rcitemark \Rspace{}.

\example{Example}  Suppose that $\G=Sp_{2n},$ and $\B=\TT\U$
is the Borel subgroup of upper triangular matrices in $\G.$
For the purposes of this discussion, we assume that $\G$
is defined with respect to the form
$$\pmatrix
&&&&&&1\\
&&&&&-1\\
&&&&.\\
&&&.\\
&&.\\
&1\\
-1\endpmatrix.$$
Then $$\TT=\left\{\pmatrix
x_1\\
&x_2\\
&&\ddots\\
&&&x_n\\
&&&&x_n^{-1}\\
&&&&&\ddots\\
&&&&&&x_2^{-1}\\
&&&&&&&x_1^{-1}
\endpmatrix\ \bigg| \ x_i\in\G_m\right\}.$$
We might denote a typical element of $\TT$ by $t(\{x_i\}).$
The root system $\P(\G,\TT)$ is of type $C_n,$
with simple roots
$\{e_i-e_{i+1}\}$ for $i=1,\dots, n-1,$ and $2e_n.$
If $\chi$ is a character of $T=\TT(F),$
then $\chi$ is of the form
$\chi=\chi_1\otimes\dots\otimes\chi_n,$
for a collection of characters $\{\chi_i\}\subset\widehat{F^\times}.$
That is,
$\chi(t(\{x_i\})=\prod\limits_i \chi_i(x_i).$

Note that if $\b=e_1-e_2,$ then $\A_\b=\{t(\{x_i\})|x_1=x_2\}.$
In this case
$$w_\b=\pmatrix 0&1\\1&0\\&&I_{2n-4}\\&&&0&1\\&&&1&0\endpmatrix.$$
So
$$\M_\b=\left\{\pmatrix
g_1\\
&x_3\\
&&\ddots\\
&&&x_3^{-1}\\
&&&& ^\tau g_1^{-1}\endpmatrix\ \bigg|\ g_1\in GL_2, x_i\in \G_m\right\}.$$
Here $^\tau g_1$ is the transpose of $g_1$ with respect to the off diagonal.
Note that
$$\N_\b=\left\{\pmatrix
\boxed{\matrix 1&x\\0&1\endmatrix}\\
&1\\
&&\ddots&&0\\
&&&1\\
&&&&\boxed{\matrix 1&-x\\0&1\endmatrix}
\endpmatrix\ \bigg|\ x\in F\right\}.$$
Set $G_1=GL_2(F),$ and let  $B_1$ be its its Borel subgroup.
Then, the analytic behavior of intertwining operator
$A(\nu,\s,w_\b)$ at $\nu=0$ detects the reducibility
of $\Ind_{B_1}^{G_1}(\chi_1\otimes\chi_2).$
The other roots of the form $e_i-e_j$ are treated in a similar fashion.

Suppose that $\a=2e_n.$
Then
$$\M_\a=\left\{\pmatrix
x_1\\
&x_2\\
&&\ddots\\
&&&h\\
&&&&\ddots\\
&&&&&x_2^{-1}\\
&&&&&&x_1^{-1}
\endpmatrix\ \bigg|\ h\in Sp_2=SL_2, x_i\in\G_m\right\}.$$
Let $H=SL_2(F),$ and $B'$ be its Borel subgroup.
Then, $A(\nu,\s,w_\a)$ can be considered as
the intertwining operator which determines the reducibility
of $\Ind_{B'}^H(\chi_n).$

As a consequence of the above computations,
we see that the Plancherel measure $\mu(\nu,\chi)$
can be computed if we understand $\mu_{G_1}(\nu_1,\chi_i\otimes\chi_j),$
and $\mu_H(\nu_2,\chi_i),$ where $\mu_{G_1}$ and $\mu_H$ have
the obvious meanings.
\endexample

Suppose $w\sigma\simeq\sigma.$
Let $T_w:V\lra V$ be an isomorphism between $w\sigma$ and $\sigma,$ i.e.,
$T_w w\sigma=\sigma T_w.$
Since
$\calA (\sigma,w)\colon \Ind_P^G(\sigma) @>>> \Ind_P^G
(w\sigma),$ we see that
$$\calA'(\sigma,w)=T_w\calA(\sigma,w):\Ind_P^G(\sigma) @>>> \Ind_P^G(\sigma).$$
We have $\calA'(\sigma,w_1 w_2)=\eta(w_1,w_2) \calA'(\sigma,w_1)\calA'
(\sigma,w_2),$ where $\eta(w_1,w_2)$ is given by $T_{w_1 w_2}=\eta(w_1,w_2)
T_{w_1}T_{w_2}.$

\proclaim{Theorem 2.15 (Harish--Chandra) Commuting Algebra Theorem}

The collection $\{\calA'(\sigma,w)|w\in W(\sigma)\}$
spans $C(\sigma).$\qed
\endproclaim

Harish-Chandra's Commuting Algebra Theorem justifies our concentration on the
theory of the operators $A(\nu,\s,w).$  We need to determine which
of the operators $\Cal A'(\s,w)$ are scalar, and which are not.
In a series of papers, beginning in the early 1960's, Kunze and Stein
computed the poles of intertwining operators for groups over $\bc$
\Lcitemark 46\Citecomma
47\Citecomma
48\Citecomma
49\Rcitemark \Rspace{}.
This work was extended by Knapp and Stein\Lspace \Lcitemark 40\Citecomma
41\Rcitemark \Rspace{}.
Finally, Knapp and Zuckerman used these results to classify the irreducible
tempered representations of Semisimple Lie Groups\Lspace \Lcitemark
42\Rcitemark \Rspace{}.
Knapp and Stein described an algorithm for determining a
basis of $C(\sigma).$
Silberger\Lspace \Lcitemark 71\Citecomma
73\Rcitemark \Rspace{} showed that this construction is valid for
$p$--adic groups. Let $\beta\in \Phi(\PPP,\A),$ and
let $\M_\beta,\ \phantom{}^*\PPP_\beta,$ and $\N_\beta$ be as before.
Let $\Delta'=\{\beta\in \Phi(\PPP,\A)|\mu_\beta(\sigma)=0\}.$
The following lemma is quite important, and its proof is non-trivial.

\proclaim{Lemma 2.16 (Knapp--Stein)} $\Delta'$ is a sub--root system of
$\Phi(\PPP,\A).$\qed
\endproclaim

Let $W'=\langle w_\beta|\beta\in \Delta'\rangle.$   The lemma guarantees that
this is well defined.  Let
$R(\s)=\{w\in W(\sigma)|w\Delta'=\Delta'\}=\{w|w\beta > 0,\ \text{for all }\
\beta
\in \Delta'\}.$  We sometimes denote $R(\s)$ by $R.$
If $\b\in\DE',$ then $w_\b\s\simeq\s,$ and
$\Ind_{^*P_\b}^{M_\b}(\s)$ is irreducible.  Therefore,
by Schur's Lemma, the normalized
operator $\Cal A'(\s,w_\b)$ is a scalar.  By the cocycle relation, we see that
$\Cal A'(\s, w)$ is scalar for every $w\in W'.$
On the other hand, if $r\in R,$ then $N_r\cap N_{w_\b}=\{I\},$
for each $\b\in W'.$  Suppose that $r=w_\a$ for some $\a\in\P(\PPP,\A).$
Then, since $w_\a\in W(\s),$ and $\a\not\in\DE',$ we see that
$\Ind_{^*P_\a}^{M_\a}(\s)$ is reducible.  Thus, by
the Commuting Algebra Theorem,
$\Cal A'(\s,w_\a)$ is non-scalar.
Now, if $r=r'w_\a,$ then, by the cocycle relation, $\Cal A'(\s,r)$
is non-scalar.

\proclaim{Theorem 2.17 (Knapp--Stein, Silberger)}
For every $\s\in\Cal E_2(M),$ $W(\sigma)=R \ltimes W',$
and $W'=\{w|\calA' (\sigma,w)\text{ is scalar}\}.$\qed
\endproclaim

So $\{\Cal A'(\s,r)\ |\ r\in R\},$ gives a basis for $C(\s).$
Note that $\eta:R\times R@>>> \bc^\times$ is a 2--cocycle of $R.$ Moreover, we
have
$$\calA' (\sigma,w_1 w_2)=\eta(w_1,w_2)\calA'(\sigma,w_1)\calA' (\sigma,w_2).$$
Thus, $C(\sigma)\simeq \bc[R]_{\eta},$ where $\bc[R]_\eta$ is the complex
group algebra of $R,$ with multiplication twisted by $\eta.$
If $\s$ is generic, then $\eta$ splits\Lspace \Lcitemark 39\Rcitemark
\Rspace{}.
For simplicity, we assume $\eta$ splits.
Let $\rho$ be an irreducible representation of $R,$ and
suppose  $\chi_\rho$ is it's
character.
Set
$$
A_\rho={1\over |R|} \dim \rho \sum_R\overline{\chi_\rho} (r) \calA'
(\sigma,r).
$$
If $\rho$ and $\rho'$ are
irreducible representations of $R,$ we have
$$\ds{A_\rho A_{\rho'}={1\over |R|^2} \dim\rho \dim\rho'
\sum\limits_{r,r'} \overline{\chi_\rho} (r) \overline{\chi_{\rho'}} (r')
\calA'(\sigma,r r')}$$
$$
\align
&={1\over |R|^2} \dim\rho \dim\rho' \sum_{w\in R} \left(\sum_{rr'=w}
\overline{\chi_\rho} (r) \overline{\chi_{\rho'}} (r')\right)\calA' (\sigma,w)\\
&={1\over |R|^2} \dim\rho \dim\rho' \sum_{w\in R} \overline{\chi_{\rho'}}(w)
\left(\sum_{r\in R} \overline{\chi_\rho} (r) \overline{\chi_{\rho'}}(r)^{-1}
\right)\calA' (\sigma,w).\\
\endalign
$$
By the Schur orthogonality relations, we see that this is
$0$ if $\rho\ncong \rho',$ and is
$A_\rho$ if $\rho\simeq \rho'.$
So the $A_\rho$ are orthogonal projections. Note that
$$
\align
\calA'(\sigma,w)A_\rho&={1\over |R|}\dim\rho\sum_r
\overline{\chi_\rho}(r)\calA' (\sigma,wr)\\
&={1\over |R|}\dim\rho\sum_r \overline{\chi_\rho} (w^{-1}rw)\calA'
(\sigma,w w^{-1}rw)\\
&={1\over |R|}\dim\rho \sum_r \chi_\rho (r) \calA' (\sigma,rw)\\
&=A_\rho \calA' (\sigma,w).
\endalign
$$
So each $\calA_\rho$ is in the center of $\bc(\sigma).$

Suppose $\Ind_P^G(\sigma)=m_1 \pi_1\oplus \ldots \oplus m_n \pi_n.$
Then $\dim C(\sigma)=m_1^2+\ldots + m_n^2.$
But $\dim \bc[R]=|R|=\sum\limits_{\rho\in\hat{R}}(\dim\rho)^2.$
Moreover, $\dim Z(\bc[R])=|\hat R|.$ Thus,
if $V_\rho=\text{Im}(\calA_\rho),$ then $V_\rho$ must be an
isotypic subspace.

\proclaim{Theorem 2.18 (Keys\Lspace \Lcitemark 39\Rcitemark \Rspace{})} Assume
that $\eta$ splits.
\item{(1)}The inequivalent components of $\Ind_P^G(\sigma)$ are parameterized
by the irreducible representations $\rho$ of $R.$
\item{(2)}$\dim \text{Hom}_G (\pi_\rho,\Ind_P^G(\sigma))=\dim\rho.$\qed
\endproclaim

\subheading{(d)  Computations} In order to compute the poles of
$A(\nu,\sigma,w),$ it is
helpful to know the following result, which follows from the fact that $PN^-$
is dense in $G,$ and
that there is a non--degenerate pairing between $V(\nu,\sigma)$ and
$V(-\nu,\sigma).$

\proclaim{Lemma 2.19 (Rallis\Lspace \Lcitemark 68\Rcitemark \Rspace{})}Suppose
$\PPP=\M\N$ is a maximal proper parabolic subgroup of $\G.$
Let ${\N}^-$ be the unipotent radical opposed to
$\N.$
Then every pole of $\nu\mapsto A(\nu,\sigma,w)$ is a pole of
$\nu \mapsto A(\nu,\sigma,w)
f(e),$ for some $f\in V(\nu,\sigma)$ with $\text{supp } f$
contained in $ PN^-$ modulo $P.$\qed
\endproclaim

\example{Example 1} Our first example is the computation of the pole of the
intertwining operator for $G=SL_2(F).$
For a complete treatment of this case, one should consult
\Lcitemark 61\Rcitemark \Rspace{}.  Let $\G=SL_2,$ and take
$\B$ to be the Borel subgroup of upper triangular matrices
of determinant one,
Let $\chi\colon F^\times @>>> \bc^x,$ be a character.
We also denote by $\chi$ the character of
$B=\B(F)$ given by $\chi\pmatrix a&b\\
0&a^{-1}\endpmatrix=\chi(a).$
Moreover, every character of $B$ is of this form, for some
$\chi\in\widehat{F^\times}.$
We have $\frak a^*_\bc/\frak z\simeq \bc,$ and we choose the isomorphism given
by
$s\mapsto |\ |^s.$  Note that $\d_{B}\left(\pmatrix
a&0\\0&a^{-1}\endpmatrix\right)=|a|^2.$
Let $I(s,\chi)=\Ind_P^G (\chi\otimes | \ |^s).$
Suppose $C$ is a compact neighborhood of 0 in $F,$ and let
$$f\left(\pmatrix 1&0\\
x&1\endpmatrix\right)=\cases1&\text{if $x\in C$}\\
0&\text{if $x\notin C$}.\endcases$$

Now extend $f$ to $B U^-,$ by $f\left(\pmatrix a&b\\
0&a^{-1}\endpmatrix \pmatrix 1&0\\
x&1\endpmatrix\right)=\chi(a)|a|^{s+1} f\pmatrix 1&0\\
x&1\endpmatrix.$
We choose  $w=\pmatrix 0&1\\
-1&0\endpmatrix$ to represent the non-trivial Weyl group element.
Then
$$
A(s,\chi,w)f(e)=\int\limits_{U}f(w^{-1} u) du.
$$

\noindent
If $u=\pmatrix 1&x\\
0&1\endpmatrix,$ then $w^{-1} u\in B U^-$ if, and only if $x\in F^\times,$
in which case
$$
w^{-1} u=\pmatrix x^{-1}&-1\\
0&x\endpmatrix \pmatrix 1&0\\
x^{-1}&1\endpmatrix.
$$
So
$$A(s,\chi,w) f(e)=\int\limits_{F^\times} \chi (x)^{-1}|x|^{-s+1}f\pmatrix
1&0\\
x^{-1}&1\endpmatrix dx=\int\limits_C \chi(x)|x|^{s-1}\ dx.$$
The poles of this last expression
come from $L$\snug-function $L(\chi,s)$ of
Tate's thesis\Lspace \Lcitemark 79\Rcitemark \Rspace{}.  We have
$L(\chi,s)=1,$ unless $\chi$ is unramified, in which case
$$
L(\chi,s)={1\over 1-\chi(\varpi)q^{-s}}.
$$
So, $A(s,\chi,w)$ has a pole at $s=0$ if and only if $L(\chi,s)$ has a pole
at $s=0,$ which occurs if and only if $\chi=1.$

Note that
$$w\chi\pmatrix a&*\\
0&a^{-1}\endpmatrix=\chi \left(w^{-1} \pmatrix a&*\\
0&a^{-1}\endpmatrix w\right)=\chi\pmatrix a^{-1}&*\\
0&a\endpmatrix=\chi(a^{-1})=\chi^{-1}(a).$$
So $w\chi\simeq \chi$ if and only if $\chi=\chi^{-1}.$  Therefore, we have
the following result.
\proclaim{Theorem 2.20 (Sally\Lspace \Lcitemark 61\Rcitemark \Rspace{})}
Let $\chi\in\widehat{F^\times},$ and extend $\chi$
to a character of $B.$  Then
$\Ind_P^G(\chi)$ is reducible if, and only if, $\chi^2=1,\ \chi\not= 1.$\qed
\endproclaim

If $\G$ is a Chevalley group\Lspace \Lcitemark 6\Citecomma
74\Rcitemark \Rspace{} and $\B$ is a minimal
parabolic
subgroup, then Winarsky,\Lspace \Lcitemark 81\Rcitemark \Rspace{}, computed the
poles of the intertwining operators
$A(\nu,\chi,w)$ for any $\chi\in\widehat B,$ and Weyl group element $w.$
Winarsky's computation reduces the computation of Plancherel measures
to the case of $SL_2,$ and then uses the techniques of\Lspace \Lcitemark
61\Rcitemark \Rspace{}.
\endexample

\example{Example 2} We next consider the case $\G=GL_n.$
This case was studied by Jacquet,\Lspace \Lcitemark 36\Rcitemark \Rspace{},
Jacquet and Godement,\Lspace \Lcitemark 24\Rcitemark \Rspace{},
Ol'{\v{s}}anski{\v{i}},
\Lcitemark 58\Rcitemark \Rspace{}, Bernstein and Zelevinsky,\Lspace \Lcitemark
4\Citecomma
5\Rcitemark \Rspace{}, and
Shahidi\Lspace \Lcitemark 65\Rcitemark \Rspace{}.
The derivation of the pole of the intertwining operator is due to
Ol'{\v{s}}anski{\v{i}}.  This result was duplicated by Bernstein and
Zelevinsky by different methods.
In\Lspace \Lcitemark 37\Rcitemark \Rspace{}, Jacquet showed that
every tempered representation of $GL_n$ is generic.
Shahidi used the local coefficient to derive an explicit formula
for the Plancherel measure.  Bushnell and Kutzko\Lspace \Lcitemark 13\Rcitemark
\Rspace{}
have recently classified the admissible dual of $GL_n,$ using the theory of
types. Their techniques also give many of the above results.

Let $\B$ be the Borel subgroup of  non-singular upper triangular matrices.
We will consider maximal parabolic subgroups. Suppose that $n=m+k,$ and that
$$\A=\left\{\pmatrix \lambda I_k&0\\0&\eta I_m\endpmatrix\bigm{|}
\lambda,\eta\in \G_m^\times\right\}.$$
Let
$$\M=Z_{\G}(\A)=\left\{\pmatrix g&0\\0&h\endpmatrix\bigm{|}\ \matrix g\in
GL_k\\h\in GL_m\endmatrix\right\}\simeq GL_k\times GL_m.$$
We take $\N=\left\{\pmatrix I_k&*\\0&I_m\endpmatrix\right\}.$
Then $\PPP=\M\N$ is a maximal parabolic subgroup of $\G,$
and every  standard maximal parabolic subgroup of $\G$ is of this form.

Let $\s=\s_1\otimes\s_2$ be an irreducible unitary supercuspidal
representation of $M.$  Let $\omega_{\s_i}$ be the central character
of $\s_i.$
Note that, if $m\neq k,$ then $W(\G,\A)=\{1\},$ so, by Bruhat's Theorem,
$\Ind_P^G(\s)$ is irreducible.  Suppose $m=k.$  Then $w=\pmatrix
0&I_k\\I_k&0\endpmatrix$ represents the unique nontrivial element of
$W(\G,\A).$
If $(g,h)\in M,$ then $w((g,h))=(h,g).$  So, $w\s\simeq\s_2\otimes\s_1,$
and $w\s\simeq\s$ if and only if $\s_1\simeq\s_2.$
Let $V$ be the space on which $\s_1$ acts, and suppose that $\s_2=\s_1.$
For  $v\in V$ and $\tilde v\in \widetilde V$ we denote
by $\ph_{v,\tilde v}$ the associated matrix coefficient,
i.e.,
$\ph_{v,\tilde v}(x)=<\tilde v,\s_1(x)v>.$
Choose a compact subset $L$ of $0$ in $M_k(F),$ and vectors
$v_1,v_2\in V.$
Let $f\left(\pmatrix I&0\\X&0\endpmatrix\right)=\xi_L(X)(v_1\otimes v_2.),$
where $\xi_L$ is the characteristic function of $L.$
A straightforward matrix computation shows that, if $n=\pmatrix
I&X\\0&I\endpmatrix,$ then $w^{-1}n\in PN^-$ if and only if $X\in GL_k.$
In this case,
$$w^{-1}n=\pmatrix -X^{-1}&I\\0&X\endpmatrix \pmatrix
I&0\\X^{-1}&I\endpmatrix.$$

Fix $\tilde v_1,\tilde v_2\in \widetilde V.$
Note that $\frak a_\bc^*/\frak z\simeq \bc,$ via $s:(g,h)\mapsto |\det
gh^{-1}|^{s/2}.$  We choose this normalization because it is the one used
in\Lspace \Lcitemark 65\Rcitemark \Rspace{} and\Lspace \Lcitemark 67\Rcitemark
\Rspace{}. This allows us to easily
describe the complementary series (see subsection (e).)
Let $dX$ be a Haar measure on $M_k(F),$ and $d^\times X$
the associated Haar measure on $GL_k(F).$

We have
$$\gather
<\tilde v_1\otimes\tilde v_2, A(s,\s,w)f(e)>=<\tilde v_1\otimes\tilde
v_2,\int\limits_N f(w^{-1}n)\, dn>\\
=\int\limits_{GL_k(F)}<\tilde v_1\otimes\tilde v_2,f\left(\pmatrix
-X^{-1}&I\\0&X\endpmatrix\pmatrix I&0\\X^{-1}&I\endpmatrix\right)>\ dX\\
=\omega_{\s_1}(-1)\int\limits_{GL_k(F)}<\tilde v_1\otimes\tilde
v_2,\s_1(X)^{-1}v_1\otimes\s_1(X)v_2>|\det X|^{-s +2k}\xi_L(X^{-1})\ dX\\
=\int\limits_{GL_k(F)}\varphi_{v_1,\tilde v_1}(X)\varphi_{v_2,\tilde
v_2}(X^{-1})\xi_L(X)|\det X|^{s} d^\times X.\tag 2.2
\endgather$$
Now, (2.2) always has a pole at s=0, for some choice of $v_i,\tilde v_i,$
and $L.$  For example, choose $v_1=v_2\neq0,$ and $\tilde v_1=\tilde v_2,$
with $<\tilde v_1,v_1>=1,$ and $L=K_0,$ an open compact subgroup with
$v_1\in V^{K_0}.$
Then
$$(2.2)\sim\int\limits_{K_0} |\det X|^{s}\ d^\times X,$$
has a pole at $s=0.$
\endexample

\example{Example 3} We compute the $R$\snug-group in a specific instance.  For
Chevalley groups, Keys,\Lspace \Lcitemark 38\Rcitemark \Rspace{}, used the
results of\Lspace \Lcitemark 81\Rcitemark \Rspace{} to
compute the $R$\snug-groups when $\PPP$ is minimal.  We examine an
example when the minimal parabolic subgroup is of $p$\snug-rank two.

Let
$$\G=Sp_4=\left\{g\in GL_4\bigg |\ ^tg\pmatrix &&&1\\&&1\\&-1\\ -1\endpmatrix
g=\pmatrix &&&1\\ &&1\\&-1\\ -1\endpmatrix\right\}.$$
Let
$$\ds{\TT=\left\{ \pmatrix x\\
&y\\
&&y^{-1}\\
&&&x^{-1}\endpmatrix\Bigg|x,y\in \G_m\right\}}.$$
We denote a typical element of $\TT$ as $t(x,y).$
$\Phi(\G,\TT)$ is of type $C_2.$
Let $\alpha=e_1-e_2,$ and $\beta=2 e_2$ be the simple roots
given by
$$
\alpha(t(x,y))=x y^{-1},\text{ and }\quad \beta(t(x,y))=y^2.
$$
Note that $W(\G,\TT)\cong \bz_2 \ltimes (\bz_2\times \bz_2),$
with generators
$$\align
w_\alpha&\colon (t(x,y))\mapsto t(y,x),\\
w_\beta&\colon t((x,y))\mapsto t(x,y^{-1}),\text{ and}\\
w_\gamma&\colon t((x,y)) \mapsto t(x^{-1},y).
\endalign
$$
Note that $w_\alpha w_\beta w_\alpha^{-1}=w_\gamma.$

Let $T=\TT(F),$ and suppose that $\chi\in\widehat T.$
Then,  for some $\chi_1,\chi_2\in\widehat{F^\times},$
we have $\chi(t(x,y))=\chi_1(x)\chi_2(y).$   Therefore,
we write $\chi=(\chi_1,\chi_2).$
Computing directly, we see that
$w_\alpha \chi=\chi $ if and only if
$\chi_1=\chi_2,$ and  $w_\beta(\chi)=\chi$ if and only if
$\chi_2^2=1.$

We have ${\A_\alpha=\left\{t(x,x)|x\in\G_m\right\}},$ and
$\M_\alpha=\left\{\pmatrix g&0\\
0&h\endpmatrix \ \bigg|\ g,h\in GL_2\right\}\cap \G.$
Since
$$
\align
&\pmatrix {}^t\! g&0\\
0&{}^t\! h\endpmatrix \pmatrix &&&1\\
&&1\\
&-1\\
-1\endpmatrix \pmatrix g&0\\0&h\endpmatrix=
\pmatrix &&&1\\
&&1\\
&-1\\
-1\endpmatrix ,\\
\intertext{we have}
& {}^t\! g \pmatrix 0&1\\
1&0\endpmatrix h=\pmatrix 0&1\\
1&0\endpmatrix \text{ and thus, }h=\pmatrix 0&1\\
1&0\endpmatrix {}^t\! g^{-1} \pmatrix 0&1\\
1&0\endpmatrix={}^\tau\! g^{-1},
\endalign
$$
where $^\tau\pmatrix a&b\\c&d\endpmatrix=\pmatrix d&b\\c&a\endpmatrix.$
So $\M_\alpha\simeq GL_2.$  Let $\G'=GL_2,$ and denote its Borel subgroup
by $\B'.$
Recall that $W(\M_\alpha,\TT)=\{1,w_\alpha\}.$
So $\chi$ ramifies in $M_\alpha$ if and only if $\chi_1=\chi_2.$
By example 2,\Lspace \Lcitemark 5\Citecomma
58\Rcitemark \Rspace{},
$\Ind_{^*B_\alpha}^{M_\alpha} (\chi)\simeq
\text{ Ind}_{B'}^{G'} (\chi_1\otimes \chi_2)$ is always irreducible.
Thus, $\mu_\alpha (\chi)=0$ if and only if $\chi_1=\chi_2.$

We next consider the simple root $\b.$
We have
$\A_\beta=\left\{ t(x,1)\big| x\in\G_m \right\}.$
Thus,
$$\M_\beta=\left\{\pmatrix x&0&0\\
0&h&0\\
0&0&x^{-1}\endpmatrix \Bigg|\matrix x\in \G_m\qquad\\
h\in GL_2 \endmatrix \right\}\cap\G.
$$
Note that we must have $^t\! h \pmatrix 0&1\\
-1&0\endpmatrix h=\pmatrix 0&1\\
-1&0\endpmatrix,$ i.e. $h\in Sp_2 =SL_2 .$
So $\M_\beta \simeq \G_m \times SL_2 .$
Let $\G''=SL_2,$ and denote its Borel subgroup by $\B''.$
We have
$W_\beta(\M_\beta,\TT)\simeq \langle w_\beta\rangle,$ and
$\Ind_{^*B_\beta}^{M_\beta} (\chi) \cong \chi_1 \otimes \text{ Ind}_{B''}^{G''}
(\chi_2)$ is reducible if
and only if $\chi_2^2=1,\ \chi_2\not= 1$ (example 1).
So $\mu_\beta (\chi)=0$ if and only if $\chi_2=1.$
Similarly, since $\gamma=2e_1=w_\alpha (\gamma),$ we have $\mu_\gamma(\chi)=0$
if and only if $\chi_1=1.$

We can now compute the $R$--groups.
The result we state is Theorem $C_n$ of\Lspace \Lcitemark 38\Rcitemark
\Rspace{}
for  $n=2.$

\proclaim{Theorem 2.21 (Keys)} Let $R$ be the $R$--group attached to $\chi.$
Then $R\simeq \bz_2^d,$ where $d$ is the number of unequal elements of
$\{\chi_1,\chi_2\}$ for which $\chi_i^2=1,$ with $\chi_i\not= 1.$
\endproclaim

\demo{Proof}
We first claim that $w_\alpha \notin R$ for any $\chi.$
Suppose $w_\alpha \in W(\chi).$  Then, $\chi_1=\chi_2,$ and thus,
$\mu_\alpha (\chi)=0,$ i.e.\ $\alpha\in \Delta'.$
Since $w_\alpha (\alpha)=-\alpha < 0,$ we have $w_\alpha\notin R,$ as claimed.
Therefore, $R\subseteq \bz_2 \times \bz_2=\langle w_\gamma,w_\beta\rangle.$
For $w\in W,$ we let $R(w)=\{\delta>0|w\delta<0\}.$
Then $R(w_\b)=\{\b\},$  $R(w_\gamma)=\{\alpha,\gamma\},$ and
$R(w_\b w_\gamma)=\{\a,\b,\gamma\}.$
Since $R(w_\b w_\gamma)=\Phi^+(\G,\TT),$
we see that $w_\b w_\gamma\in R$ implies $\DE'=\emptyset.$
However, $w_\b w_\gamma\in W(\chi)$ also implies that
$\chi_1^2=\chi_2^2=1,$
and thus,
$w_\gamma, w_\b\in W(\chi).$ Therefore, we see that
$w_\b w_\gamma\in R$ implies $w_\b$ and $w_\gamma\in R.$
$w_\gamma,w_\b\in R.$
Since $R(w_\b)=\{\b\},$ we have
$w_\beta \in R $ if and only if  $w_\beta \in W(\chi),$
and $\mu_\beta (\chi)\not= 0,$ i.e.\ $\chi_2^2=1,\ \chi_2\not= 1.$
If $\chi_1=\chi_2,$ and $w_\gamma \in W (\chi),$ we have
$\alpha\in \Delta',$ and thus $w_\gamma \notin R.$
Therefore, $w_\gamma \in R$ if and only if $\chi_1^2=1,\ \chi_1\not= 1,$ and
$\chi_1\not= \chi_2.$ This completes the proof\qed
\enddemo
\endexample

Knapp and Zuckerman,\Lspace \Lcitemark 43\Rcitemark \Rspace{}, were first to
find
an example of a non-abelian $R$\snug-group, showing that
sometimes $\Ind_P^G(\s)$ has components which appear with multiplicity larger
than 1.  Keys found many more examples in\Lspace \Lcitemark 38\Citecomma
39\Rcitemark \Rspace{}.
For $G=Sp_{2n}(F),$ or $SO_n(F),$ Goldberg computed the $R$\snug-groups
for all parabolic subgroups in\Lspace \Lcitemark 26\Citecomma
27\Rcitemark \Rspace{}.
For $SL_n$ the possible $R$\snug-groups are computed in\Lspace \Lcitemark
29\Rcitemark \Rspace{},
which builds upon the case of the minimal parabolic, where the $R$\snug-groups
were known from\Lspace \Lcitemark 20\Citecomma
21\Citecomma
38\Rcitemark \Rspace{}.  For the quasi-split unitary groups
$U_n,$ Goldberg computed the $R$\snug-groups in\Lspace \Lcitemark 25\Rcitemark
\Rspace{}.

\subheading{(e) Elliptic Representations}
The collection of irreducible elliptic
tempered representations played a key role
in determining the tempered spectrum of real reductive groups\Lspace \Lcitemark
42\Rcitemark \Rspace{}.  While their analog in the $p$\snug-adic case are
important,
understanding their characters will not be enough to determine the characters
of all tempered representations.  We describe what is known about such
 representations.

A regular element\Lspace \Lcitemark 57\Rcitemark \Rspace{} of $G$ is said to be
{\bf elliptic} if its
centralizer is compact modulo the center if $G.$  We write $G^e$ for the
set of regular elliptic elements of $G.$  An irreducible admissible
representation
$\pi$ of $G$ is said to be {\bf elliptic} if $\chi_\pi$ is not identically
zero on $G^e.$  Here, $\chi_\pi$ is the distribution character of $\pi$
\Lcitemark 11\Rcitemark \Rspace{}. The following result exhibits the importance
of elliptic representations.

\proclaim{Theorem 2.22 (Knapp-Zuckerman\Lspace \Lcitemark 42\Rcitemark
\Rspace{})}
Suppose that $\G$ is a connected reductive algebraic group defined
over $\R,$ and set $G=\G(\R).$
Let $\pi$ be an irreducible tempered representation of $G.$
If $\pi$ is not elliptic, then there is some proper
parabolic subgroup $\PPP=\M\N$ of $G,$ and an irreducible elliptic tempered
representation $\s$ of $M$ for which
$\pi=\Ind_P^G(\s).$\qed
\endproclaim

We now consider the case where $F$ is  a $p$\snug-adic field
of characteristic zero.
Arthur,\Lspace \Lcitemark 3\Rcitemark \Rspace{}, has given a necessary and
sufficient condition,
in terms of $R$\snug-groups,
for a tempered representation to be elliptic.  For simplicity, we assume that
the $R$\snug-group in question is abelian, and the cocycle $\eta$ is a
coboundary. (Arthur makes no such assumptions, and his result
becomes slightly more technical.)
For $H\in\frak a,$ we let $w\cdot H$ denote the image of
$H$ under
the action of $W(\G,\A)$ on $\frak a.$  Let
$$\frak a_w=\{H\in\frak a|w\cdot H=H\}.$$
Suppose $\s\in\Cal E_2(M),$ and suppose that $R=R(\s)$
is the $R$\snug-group attached to $\Ind_P^G(\s).$   Then we let
$\frak a_R=\bigcap\limits_{w\in R} a_w.$

\proclaim{Theorem 2.23 (Arthur\Lspace \Lcitemark 3\Rcitemark \Rspace{})}
Suppose that $\s\in\Cal E_2(M),$
and that $R=R(\s)$ is abelian.  We further assume that the $2$\snug-cocycle
$\eta$ attacher to $\s$ and $R$ splits.  Then the following are equivalent:
\roster
\item"(a)" $\Ind_P^G(\s)$ has an elliptic component;
\item"(b)" Every component of $\Ind_P^G(\s)$ is elliptic;
\item"(c)" There is a $w\in R,$ with $\frak a_w=\frak z.$\qed
\endroster
\endproclaim

In view of the theorem of Knapp and Zuckerman, it is reasonable to ask which
irreducible tempered representations of $G$ are irreducibly induced from
elliptic representations of proper Levi subgroups.  In\Lspace \Lcitemark
33\Rcitemark \Rspace{} Herb
gives a description of such representations, within the constraints of the
previous theorem.

\proclaim{Theorem 2.24 (Herb\Lspace \Lcitemark 33\Rcitemark \Rspace{})}  Let
$\s,\ R,$ and $\eta$
be as above.  Let $\pi$ be a subrepresentation of
$\Ind_P^G(\s).$  Then there is a parabolic subgroup $\PPP'=\M'\N',$
and an irreducible elliptic  $\tau\in\Cal E_t(M')$ with
$\pi=\Ind_{P'}^G(\tau),$
if and only if there is a $w\in R$ with $\frak a_w=\frak a_R.$\qed
\endproclaim

Herb was able to use these results to describe the elliptic tempered
representations of $Sp_{2n}$ and $SO_n$\Lspace \Lcitemark 33\Rcitemark
\Rspace{}. For $Sp_{2n}(F)$
and $SO_{2n+1}(F)$ the analog of  Theorem 2.22 holds.  That is,
an irreducible tempered representation is either elliptic, or is irreducibly
induced from an elliptic representation of a proper parabolic subgroup.
For $SO_{2n}(F),$ this statement is false.  We give a simple example
due to Herb\Lspace \Lcitemark 33\Rcitemark \Rspace{}.
Let $\G=SO_6.$  We define $\G$ with respect to the form
$\pmatrix
&&&&&1\\
&&&&1\\
&&&1\\
&&1\\
&1\\
1\endpmatrix.$
Then
$$\TT=\left\{
\pmatrix
x_1\\
&x_2\\=
&&x_3\\
&&&x_3^{-1}\\
&&&&x_2^{-1}\\
&&&&&x_1^{-1}\\
\endpmatrix\bigg|\ x_i\in\G_m\right\}.$$
The root system
$\P(\G,\TT)$ is of type $D_3,$ and the Weyl group is isomorphic
to $S_3\ltimes (\bz_2\times\bz_2).$  Here $S_3$ acts on the indices of
the $x_i$\snug's.  The subgroup $\{1\}\ltimes(\bz_2\times\bz_2)$
is given by
$\{1,c_1c_2,c_2c_3,c_1c_3\},$
where
$$\gather
c_1=\pmatrix
&&1\\
&I_4\\
1\endpmatrix,\\
c_2=\pmatrix
1\\
&0&&&1\\
&&1\\
&&&1\\
&1&&&0\\
&&&&&1\endpmatrix,\\
\intertext{and}
c_3=\pmatrix
I_2\\
&0&1\\
&1&0\\
&&&I_2\endpmatrix.
\endgather$$

Let $\chi\in\hat T.$  Then for some
$\chi_1,\chi_2,\chi_3\in\widehat{F^\times},$
we have $\chi=\chi_1\otimes\chi_2\otimes\chi_3.$
We need the following lemma from\Lspace \Lcitemark 38\Rcitemark \Rspace{}.
This result follows
from example (2),
and its proof is similar to the proof that $w_\a\not\in R$ for
any $\chi$ in example 3.

\proclaim{Lemma 2.25 (Keys)}  If $w=sc\in R,$ with $s\in S_3$ and $c\in
\bz_2\times\bz_2,$
then $s=1.$\qed
\endproclaim

Therefore, $R\subset \bz_2\times\bz_2$ is abelian.
Let $\chi=\chi_1\otimes\chi_2\otimes\chi_3\in\hat T.$
Then $c_ic_j\in W(\chi)$ if and only if $\chi_i^2=1,$
and $\chi_j^2=1.$
Let $d$ be the number of unequal elements of
$\{\chi_1,\chi_2,\chi_3\}$ which satisfy $\chi_i^2=1.$
If $d=0,$ then $W(\chi)=\{1\},$ so $R=\{1\}.$
Otherwise, $\R\simeq \bz_2^{d-1}.$  In particular,
if $\chi_i^2=1,$ for $i=1,2,3,$ and $\chi_i\neq \chi_j,$
for $i\neq j,$ then $R=\bz_2\times\bz_2.$  Suppose that this is the case.
Then $\Ind_B^G(\chi)$ has four inequivalent components.
Note that $\frak a=\left\{\text{diag}\{a,b,c,-c,-b,-a\}|a,b,c\in\R\right\}.$
Furthermore,
$\frak z=\{0\}.$
We denote a typical element of $\frak a$ by $t(a,b,c).$
Direct computation shows that $\frak a_{c_1c_2}=\{t(0,0,c)|c\in\R\},$
$\frak a_{c_1c_3}=\{t(0,b,0)|b\in\R\},$ and  $\frak
a_{c_2c_3}=\{t(a,0,0)|a\in\R\}.$  Therefore, $\frak a_w\neq\frak z,$
for all $w\in R.$  Consequently,  $\Ind_B^G(\chi)$ does not have elliptic
constituents.  Moreover, since $\frak a_R=\{0\},$  these components
cannot be irreducibly induced from an elliptic
tempered representation of a proper parabolic subgroup.

In fact, what happens is that whenever $\chi$ is induced to
a rank one parabolic subgroup, P=MN, the induced representation
breaks into two elliptic components.
Since there are four components of $\Ind_B^G(\chi),$
each of the components of $\Ind_{B}^M(\chi)$ induces to $G$
reducibly.  Since $SO_6$ is locally homeomorphic to $SL_4,$
such a construction is valid for the Borel subgroup of $SL_4$ as well.\qed

Herb classified  all the irreducible tempered representations of
$SO_{2n}(F)$ which are non-elliptic, and cannot be irreducibly induced from
elliptic representations.    Goldberg
classified all the possible $R$\snug-groups for $SL_n,$
and, motivated by\Lspace \Lcitemark 33\Rcitemark \Rspace{}, determined all
irreducible tempered representations which are
non-elliptic, and cannot be irreducibly induced from
elliptic representations\Lspace \Lcitemark 29\Rcitemark \Rspace{}.  For $U_n$
Goldberg showed
that such representations do not exist\Lspace \Lcitemark 25\Rcitemark
\Rspace{}.

\subheading{(f) Complementary series}
So far we have discussed the structure of $\Ind_P^G(\s),$
with $\s\in\Cal E_2(M).$  However, one would like to know the structure
of the representations $I(\nu,\s).$  For instance, one would
like to know when $I(\nu,\s)$ is irreducible and unitarizable.
Such representations are said to be in the  complementary series.  Also,
one would like to know the points of reducibility for $I(\nu,\s),$
and what properties its subrepresentations and subquotients
possess.  We outline now what is known.
For $G=GL_n,$  and $P$ maximal,
Bernstein and Zelevinsky\Lspace \Lcitemark 5\Citecomma
82\Rcitemark \Rspace{}  showed
that $\Ind_P^G(\s_1\otimes\s_2)$ is reducible if and only
if $\s_2\simeq\s_1\otimes |\ |^{\pm 1}.$
Bernstein and Zelevinsky showed that
at the points of reducibility, there is a unique non-supercuspidal
discrete series representations.  These representations are often
referred to as {\bf special} representations.
Shahidi\Lspace \Lcitemark 65\Rcitemark \Rspace{}
proved that the local coefficients satisfy the relation
$$C_\psi(s,\s_1\otimes\widetilde\s_2)C_\psi(1-s,\widetilde\s_1\otimes\s_2)=\omega_{\s_1}^k\omega_{\s_2}^k(-1).\tag 2.3$$
If $\s_1=\s_2$ is unitary, then $A(s,\s,w)$ has a pole at $s=0,$
so $\mu(s,\s_1\otimes\s_1)$ has a zero at $s=0.$
Thus, by Shahidi's theorem on local coefficients, and (2.3),
we see that $\mu(s,\s_1\otimes\s_1)$ must have a pole at $s=1.$
Thus, $I(1,\s_1\times\s_1)$ is reducible, duplicating
the results of Bernstein and Zelevinsky for maximal parabolic subgroups.

Shahidi was able to prove a general result about the complementary series
when $P=MN$ is maximal.  Suppose that $^LP=\ ^LM\phantom{}^LN$
is the standard parabolic subgroup of $^LG$ with Levi subgroup $^LM$
\Lcitemark 7\Rcitemark \Rspace{}.
Let $^L\frak n$ be the real Lie algebra of $^LN.$
Then $^LM$ acts on $^LN$ by the adjoint representation.  We denote
this representation by $r.$
There is a particular way of ordering the components of
$r,$ as described in\Lspace \Lcitemark 67\Rcitemark \Rspace{}, and we write
$r=r_1\oplus\dots\oplus r_m,$ accordingly.
We choose the isomorphism $\frak a^*_\bc/\frak z\simeq\bc$ as
in\Lspace \Lcitemark 67\Rcitemark \Rspace{}.  Let $\s$ be an irreducible
unitary generic supercuspidal representation of $M.$
We let $P_{\s,r_i}$ be the unique
polynomial satisfying $P_{\s,i}(0)=1,$ and $P_{\s,i}(q^{-s})$ is the
numerator of $\gamma(s,\s,r_i,\psi_F).$
We set $L(s,\s,r_i)=P_{\s,i}(q^{-s})^{-1}.$

We now describe a convenient parameterization of $\frak a_\bC^*.$
Let $\a$ be the unique simple root in $\N,$
and let $\rho_{\PPP}$ be half the sum of the positive roots in $\N.$
Let
$$<\rho_{\PPP},\a>=2(\rho_{\PPP},\a)/(\a,a),$$
where $(\, ,)$ is the standard euclidean inner product
on $\P(\G,\TT).$   We define an element $\wt \a$
of $\frak a_\bC^*$ by
$\wt\a=<\rho_{\PPP},\a>^{-1}\rho_{\PPP}.$  We
let $I(s,\s)=I(s\wt\a,\s).$

\proclaim{Theorem 2.26 (Shahidi\Lspace \Lcitemark 67\Rcitemark \Rspace{})}
Let $\PPP=\M\N$ be a maximal parabolic subgroup of $\G.$
Assume that $W(\G,\A)\neq\{1\},$ and let $w_0$ represent
the unique non-trivial Weyl group element.
Suppose $\s\in\phantom{ }^\circ\Cal E(M)$ is generic.
\roster
\item For $3\leq i\leq m,$  we have $L(s,\s,r_i)=1.$
\item The following are equivalent
\indent
\item"a)" $s=0$ is a pole of $A(s,\s,w_0).$
\indent
\item"b)" $P_{\s,i}(1)=0$ for either $i=1,$ or $2,$
and only for one of them.
\indent
\item"c)" $w_0\s\simeq\s,$ and $\Ind_P^G(\s)$ is irreducible.
\item"(3)"  Suppose (2a) is satisfied.  Moreover assume that in (2b),
$P_{\s,i}(1)=0.$ Then:
\indent
\item"a)" For $0<s<1/i,$ The representation $I(s,\s)$ is irreducible
and unitarizable, i.e., is in the complementary series.
\indent
\item"b)" The representation $I(1/i,\s)$ is reducible, with a unique
generic special subrepresentation.  Its Langlands quotient is non-generic,
unitarizable, and non-tempered.
\indent
\item"c)" For $s>1/i,$ the representation $I(s,\s)$
is always irreducible and never unitarizable.
\item"(4)"  If $\Ind_P^G(\s)$ is reducible, then
for all $s>0,$ the representation $I(s,\s)$ is irreducible, and
is never unitarizable.\qed
\endroster
\endproclaim

Shahidi's result points
out the power of Langlands's conjectures, because
of the connection they give between number theory and harmonic analysis.
Tadic has been able to derive many
results, complementary to Shahidi's,
for classical groups via the theory of Jacquet modules
\Lcitemark 75\Citecomma
76\Citecomma
77\Rcitemark \Rspace{}.

\subheading{\S 3 Endoscopy}
The theory of twisted endoscopy has proved very useful
in the computation of the poles of intertwining operators.
The theory itself is quite technical, and therefore, beyond
the scope of these lectures.  We will attempt to illustrate its
power by working an example.  For more details on the general theory
one should consult\Lspace \Lcitemark 44\Citecomma
45\Citecomma
55\Rcitemark \Rspace{}.

Let $n\geq1,$ and $\G=Sp_{2n}, SO_{2n},$ or $SO_{2n+1}.$
In each case, $\G$ has a maximal parabolic subgroup
(called the Siegel parabolic subgroup) with $\M\simeq GL_n.$
More precisely, if $g\in GL_n(F),$
then
$$m(g)=\cases\pmatrix g&0\\0&^tg^{-1}\endpmatrix&\text{ if }\G=Sp_{2n}\text{ or
}SO_{2n},\\
\pmatrix g&0&0\\0&1&0\\0&0&^tg^{-1}\endpmatrix&\text{ if } \G=SO_{2n+1},
\endcases$$
is the corresponding element of $M=\M(F).$
Let $(\s,V)\in\Cal E_2(GL_n(F)).$
For $s\in\bc,$ we set
$I(s,\s)=\Ind_P^G(\s\otimes |\ |^s)=$
$$\{f:G\lra V| f(m(g_0)ng)=\s(g_0)|\det g_0|^{s+\delta}f(g)\},$$
where $\delta = (n+1)/2,\ (n-1)/2,$ or $n/2$ respectively.
In each case we take $w_0$ to be the longest element of the
Weyl group.

Consider the standard representation $\rho_n:GL_n(\bc)\lra GL_n(\bc).$
This is just the natural action of $GL_n(\bc)$ on $\bc^n.$
We denote the exterior square of $\rho_n$ by $\wedge^2\rho_n.$
If $v,u\in \bc^n,$ then $\wedge^2\rho_n(g)(v\wedge u)=gv\wedge gu.$
Let $\G=Sp_{2n}.$  Then $^LG^0=SO_{2n+1}(\bc),$
and the adjoint representation $r$ of $^LM$ on $^L\frak n$
is isomorphic to $\rho_n\oplus\wedge^2\rho_n.$  By
\Lcitemark 24\Rcitemark \Rspace{}, $L(s,\s,\rho_n)=1.$  So we must
examine $L(s,\s,\wedge^2\rho_n).$  More specifically,
we need to find the polynomial $P(t)$ such that
$P(0)=1$ and $P(q^{-2s})A(s,\s,w_0)$ is holomorphic and non-zero.
Then, by\Lspace \Lcitemark 67\Rcitemark \Rspace{},
$$L(s,\s,\wedge^2\rho_n)=P(q^{-s})^{-1}.$$
If $\G=SO_{2n},$
then $r=\wedge^2\rho_n,$ and if $\G=SO_{2n+1},$ then
$r=\text{Sym}^2\rho_n,$ the symmetric square of $\rho_n.$

Suppose that $\s$ is supercuspidal, and
$\varphi:W_F\lra GL_n(\bc)=\ ^LM^0,$
is the (conjectural) Langlands parameter for $\s.$
Let $N=n(n-1)/2,$ and
consider $\wedge^2\rho_n\circ\varphi:W_F\lra GL_N(\bc).$
We must have $L(s,\s,\wedge^2\rho_n)=L(s,\wedge^2\rho_n\circ\varphi),$
where this last object is the Artin $L$\snug-function.
Now, $P(1)=0$ if and only is $L(s,\wedge^2\rho_n\circ\varphi)$
has a pole at $s=0,$ which can only occur if the trivial
representation appears in $\wedge^2\rho_n\circ\varphi.$
Thus, $\text{Im }\varphi$ must fix a vector in the alternating
space $\wedge^2\bc^n.$  By duality, there must be some $B\in
(\wedge^2\bc^n)^*,$
which is fixed by $\text{Im }\varphi,$ i.e. the image of $\varphi$
must fix a skew symmetric form in $n$ variables.  Since $\s$ is supercuspidal,
$\text{Im }\varphi$ must be irreducible\Lspace \Lcitemark 7\Rcitemark
\Rspace{},
and hence $B$ must be non-degenerate.  Thus, $n$ must be even,
and $\varphi$ should factor through $Sp_n(\bc).$  Dually, one
expects $\s$ to ``come from'' $H=SO_{n+1}(F),$ since
$^LH^0=Sp_n(\bc).$

Note that we must have the following result.

\proclaim{Proposition 3.1 (Shahidi)}
If $n$ is odd, then $L(s,\s,\wedge^2\rho_n)=1.$
\endproclaim
\demo{Proof} Consider $\G=SO_{2n}$ and $\M=GL_n.$  Then
the Weyl group $W(G,A)=\{1\},$ so, by\Lspace \Lcitemark 67\Rcitemark \Rspace{},
$P(t)=1.$\qed
\enddemo

If $\G=Sp_{2n},$ then $W(G,A)\simeq \Bbb Z_2,$ and $w_0(m(g_0))=m(^tg_0^{-1}).$
Thus, $w_0\s\simeq\s$ if and only if $\s\simeq\tilde\s$\Lspace \Lcitemark
4\LIcitemark{}, \S 7\RIcitemark \Rcitemark \Rspace{}.
When the residual characteristic of $F$ is 2, then
such representations always exist.

\proclaim{Corollary 3.2 (Shahidi)}
If $n$ is odd and $\s\simeq\tilde\s,$ then $\Ind_P^G(\s)$ is reducible,
and for all $s>0,$ the representation $I(s,\s)$ is irreducible.\qed
\endproclaim

So, we need to concentrate on the case where $n$ is even.
Suppose $\G=SO_{2n},$ with $n$ even.  Then
$$N=\{\pmatrix I&X\\0&I\endpmatrix\ |\ ^tX=-X\}.$$
We have $w_0=\pmatrix 0&I\\I&0\endpmatrix.$
If $n=\pmatrix I&X\\0&I\endpmatrix\in N,$ then
$w_0^{-1}n\in PN^-$ if and only if $X\in GL_n(F).$
In this case
$$w_0^{-1}n=\pmatrix -X^{-1}&I\\0&X\endpmatrix\pmatrix
I&0\\X^{-1}&I\endpmatrix.$$
Let $L$ be compact and $h(\bar n)=\xi_L(X)v,$ for some $v\in V.$
Taking $\tilde v\in \widetilde V,$ we have
$$<\tilde v, A(s,\s,w_0)f(e)>=\int\limits_{^tX=-X\atop\det X\neq 0}<\tilde
v,\s(-X)^{-1}v>|\det X|^{-s+(n-1)/2}\xi_L(X^{-1}) dX,$$
which we rewrite as
$$\omega_\s(-1)\int\limits_{^tX=-X\atop \det X\neq 0}\ph_{v,\tilde v}(X)|\det
X|^s\xi_L(X)\ d^\times X.\tag 3.1$$

We can choose
$f\in C_c^\infty(GL_n(F))$ so that
$$\ph_{v,\tilde v}(g)=\int\limits_{Z_n(F)}f(zg)\omega_\s^{-1}(z)\ dz,$$
where $Z_n(F)$ is the center of $GL_n(F).$

Let $\gamma_0=\pmatrix 0&I\\-I&0\endpmatrix.$ Since all
symplectic forms are equivalent, any $X$ with $^tX=-X$
and $\det X\neq 0,$ is of the form $^tg\gamma_0g,$
for some $g\in GL_n(F).$
For technical reasons
we take $\t^*(g)=w^{-1}\ ^tg^{-1}w,$
with $$w=\pmatrix
&&&&&&1\\
&&&&&-1\\
&&&&.\\
&&&.\\
&&.\\
&1\\
-1\endpmatrix.$$
For $\gamma\in GL_n(F),$ we let $\G_{\t^*,\gamma}=\{g\in GL_n\ |g^{-1}\gamma
\t^*(g)=\gamma\}.$
Now  (3.1) is proportional to
$$\int\limits_{\G_{\t^*,\gamma}(F)^0\backslash GL_n(F)}\int\limits_{F^\times}
f(g^{-1}z\gamma_0\t^*(g))|\det
(g^{-1}z\gamma_0\t^*(g))|^s\xi_L(g^{-1}z\gamma_0\t^*(g))\ dz\ d\dot g.\tag
3.2$$

\def\P{\Phi}
We now define some terms.
For $\gamma\in GL_n(F),$ and $f\in C_c^\infty(GL_n(F)),$ let
$$\P_{\t^*}(\gamma,f)=\int\limits_{\G^0_{\t^*,\gamma}\backslash
GL_n(F)}f(g^{-1}\gamma\t^*(g))\ d\dot g.$$
We say that $\gamma$ is $\t^*$\snug-semisimple
if $(g,\t)$ is semisimple in the
disconnected algebraic group $GL_n\rtimes <\t^*>.$  We say
that such a
$\gamma$ is strongly $\t^*$-regular if $\G_{\t^*,\gamma}$ is abelian.
We say that $\gamma$ and $\gamma'$ are $\t^*$-conjugate if, for some
$g\in GL_n(F),\ \gamma'=g^{-1}\gamma\t^*(g).$  Finally,
$\gamma$ and $\gamma'$ are stably $\t^*$\snug-conjugate
if, for some $g\in GL_n(\bar F),$
$\gamma'=g^{-1}\gamma \t^*(g).$
We write $\gamma\sim\gamma'$ for stable $\t^*$\snug-conjugacy.
If $\gamma\sim\gamma'$ are
strongly $\t^*$\snug-regular,
then
$\G_{\t^*,\gamma}^\circ$ and $\G_{\t^*,\gamma'}^\circ$
are inner.  We can therefore transfer measures.
Set
$$\P_{\t^*}^{st}(\gamma,f)=\sum_{\gamma'\sim\gamma}\P_{\t^*}(\gamma',f).$$
This is called the $\t^*$-stable twisted orbital integral of $f$ at $\gamma.$

Let $\bold H=SO_{n+1},$ and
let $\text{\bf T}=\text{\bf T}_{\bold H}$ be a Cartan subgroup
of $\text{\bf H}$ defined over $F.$
Suppose that $\text{\bf T}'$ is a $\t^*$-stable Cartan of $GL_n,$
defined over $F.$  Let
$$\text{\bf T}'_{\t^*}=\text{\bf T}'/(1-\t^*)\text{\bf T}'.$$
There exists  $\text{\bf T}=\text{\bf T}_{\text{\bf H}},$ such that
there is an isomorphism between
$\text{\bf T}$ and $\text{\bf T}'_{\t^*}$ defined over $F.$
Thus, we get a one-to-one map $\Cal A$ between semisimple conjugacy classes
of $\text{\bf H}(\bar F),$ and $\t^*$-semisimple $\t^*$-conjugacy
classes of $GL_n(\bar F).$

\definition{Definition 3.3}  We say that $\delta\in H$ is a {\bf norm}
of $\gamma\in GL_n(F)$ if the $GL_n(\bar F)-\t^*$\snug-conjugacy
class of $\gamma$ is the image of $\delta$ under $\Cal A.$
We write $\delta =\Cal N\gamma.$
\enddefinition
If $\gamma$  is strongly $\t^*$-regular, then $\delta=\Cal N\gamma$
is strongly regular.
Suppose $\psi\in C_c^\infty(H).$  Let
$$\P(\delta,\psi)=\int\limits_{\text{\bf H}_{\delta}^0(F)\backslash
H}\psi(h^{-1}\delta h)\ d\dot h,$$
and
$$\P^{st}(\delta,\psi)=\sum_{\delta'\sim\delta}\P(\delta',\psi),$$
for $\delta$ strongly regular.

\example{Assumption 3.4}
For every $f\in C_c^\infty(GL_n(F)),$ there exists a $f^H\in C_c^\infty(H),$
so that
$$\P_{\t^*}^{st}(\gamma,f)=\P^{st}(N\gamma,f^H),$$
for every strongly $\t$\snug-regular $\gamma\in GL_n(F),$ and
$\P^{st}(\delta,f^H)=0,$ if $\delta$ is not a norm.
\endexample

\proclaim{Proposition 3.5}  Up to a non-zero constant,
$$\int\limits_{Sp_n(F)\backslash GL_n(F)} f(^tgw^{-1}gw)\ d\dot g= f^H(e),$$
for $f\in C_c^\infty(GL_n(F)).$\qed
\endproclaim

\definition{Definition 3.6} An irreducible supercuspidal
representation $\s$ of $GL_n(F)$ is said to ``come from''
$SO_{n+1}(F)$ if, $n$ is even, and $f^H(e)\neq 0,$ for
some $f\in C_c^\infty(GL_n(F))$ defining
a matrix coefficient of $\s.$
\enddefinition

\definition{Definition 3.7} We say $\s\simeq\tilde \s$ comes from $SO_n^*(F)$
if $n$ is even, and $\s$ does not come from $SO_{n+1}(F).$
We say $\s$ comes from $Sp_{n-1}(F)$ if $n$ is odd.
\enddefinition

\proclaim{Theorem 3.8}  The residue of $A(s,\s,w_0)$ at $s=0,$ given by (3.2),
is proportional to
$$\int\limits_{Sp_n(F)\backslash GL_n(F)} f(^tgw^{-1}gw)\ d\dot g.$$
Therefore, $A(s,\s,w_0)$ has a pole at $s=0$ if
and only if, for some choice of $f,$ defining a matrix coefficient
of $\s,$
$$\int\limits_{Sp_n(F)\backslash GL_n(F)} f(^tgw^{-1}gw)\ d\dot g\neq 0.$$
It follows that $\s\simeq\tilde \s, $ and $\omega_\s=1.$\qed
\endproclaim

\proclaim{Theorem 3.9 (Shahidi\Lspace \Lcitemark 68\Rcitemark \Rspace{})}
Suppose that $\s$ is an irreducible supercuspidal representation
of $GL_n(F),$ with $\s\simeq\tilde \s.$

\item{(a)} If $\G=SO_{2n+1},$  then $\Ind_P^G(\s)$ is irreducible if
and only if $\s$ comes from $SO_n^*(F),$ or $Sp_{n-1}(F).$

\item{(b)}  If $\G=Sp_{2n},$ then $\Ind_P^G(\s)$ is irreducible
if and only if $\s$ comes from $SO_{n+1}(F).$

\item{(c)} If $\G=SO_{2n},$ then $\Ind_P^G(\s)$ is irreducible if and
only if $\s$ comes from $SO_{n+1}(F)$ or $Sp_{n-1}(F).$\qed
\endproclaim

The content of Shahidi's Theorem should not be overstated.  It
says that if the correspondence $f\mapsto f^H$ exists, with the desired
properties, then $f^H(e)\neq 0$
if and only if
$$\int\limits_{Sp_n(F)\backslash GL_n(F)} f(^tgw^{-1}gw)\ d\dot g\neq0.$$
In particular, it does not guarantee the existence of the
correspondence $f\mapsto f^H.$

In\Lspace \Lcitemark 28\Rcitemark \Rspace{},
Goldberg used Shahidi's method,
and  the explicit matching between $GL_2(E)$ and $U(2),$\Lspace \Lcitemark
60\Rcitemark \Rspace{},
to describe the poles of $A(s,\s,w)$ when
$\G=U(2,2).$
He was able to show that the poles of the intertwining operator
distinguish the images of the two base change maps described
in\Lspace \Lcitemark 60\Rcitemark \Rspace{}.
In \Lspace \Lcitemark 30\Rcitemark \Rspace{} Goldberg carried out a similar
computation for $G=U(n,n),$
and gave a relation between reducibility and lifting.  Shahidi,
\Lcitemark 63\Rcitemark \Rspace{}, has studied the residues of the intertwining
operators for more maximal parabolic subgroups of $SO_{2n}(F).$
The case of a general classical group is the subject of a work in progress
by Goldberg and Shahidi.
\Refs

\message{REFERENCE LIST}

\bgroup\Resetstrings%
\def\Ecnt{4}\def\acnt{0}%
\def\Ftest{ }\def\Fstr{1}%
\def\Atest{ }\def\Astr{J\Initper  Arthur}%
\def\Ttest{ }\def\Tstr{On some problems suggested by the trace formula}%
\def\Btest{ }\def\Bstr{Lie Group Representations II}%
\def\Etest{ }\def\Estr{R\Initper  Herb%
  \Ecomma S\Initper  Kudla%
  \Ecomma R\Initper  Lipsman%
  \Eandd J\Initper  Rosenberg}%
\def\Itest{ }\def\Istr{Springer-Verlag}%
\def\Ctest{ }\def\Cstr{New York--Heidelberg--Berlin}%
\def\Stest{ }\def\Sstr{Lecture Notes in Math.}%
\def\Ntest{ }\def\Nstr{1041}%
\def\Dtest{ }\def\Dstr{1983}%
\def\Ptest{ }\def\Pstr{1--49}%
\Refformat\egroup%

\bgroup\Resetstrings%
\def\Ecnt{0}\def\acnt{0}%
\def\Ftest{ }\def\Fstr{2}%
\def\Atest{ }\def\Astr{J\Initper  Arthur}%
\def\Ttest{ }\def\Tstr{Unipotent automorphic representations: conjectures}%
\def\Jtest{ }\def\Jstr{Societ\'e Math\'ematique de France, Ast\'erisque}%
\def\Vtest{ }\def\Vstr{171-172}%
\def\Dtest{ }\def\Dstr{1989}%
\def\Ptest{ }\def\Pstr{13--71}%
\def\Astr{\Underlinemark}%
\Refformat\egroup%

\bgroup\Resetstrings%
\def\Ecnt{0}\def\acnt{0}%
\def\Ftest{ }\def\Fstr{3}%
\def\Atest{ }\def\Astr{J\Initper  Arthur}%
\def\Ttest{ }\def\Tstr{On elliptic tempered characters}%
\def\Jtest{ }\def\Jstr{Acta Math.}%
\def\Vtest{ }\def\Vstr{171}%
\def\Dtest{ }\def\Dstr{1993}%
\def\Ptest{ }\def\Pstr{73--138}%
\def\Astr{\Underlinemark}%
\Refformat\egroup%

\bgroup\Resetstrings%
\def\Ecnt{0}\def\acnt{0}%
\def\Ftest{ }\def\Fstr{4}%
\def\Atest{ }\def\Astr{I\Initper \Initgap N\Initper  Bernstein%
  \Aand A\Initper \Initgap V\Initper  Zelevinsky}%
\def\Ttest{ }\def\Tstr{Representations of the group $GL(n,F)$ where $F$ is a
local non-archimedean local field}%
\def\Jtest{ }\def\Jstr{Russian Math. Surveys}%
\def\Vtest{ }\def\Vstr{33}%
\def\Ptest{ }\def\Pstr{1--68}%
\def\Dtest{ }\def\Dstr{1976}%
\Refformat\egroup%

\bgroup\Resetstrings%
\def\Ecnt{0}\def\acnt{0}%
\def\Ftest{ }\def\Fstr{5}%
\def\Atest{ }\def\Astr{I\Initper \Initgap N\Initper  Bernstein%
  \Aand A\Initper \Initgap V\Initper  Zelevinsky}%
\def\Ttest{ }\def\Tstr{Induced representations of reductive $p$\snug-adic
groups. I}%
\def\Jtest{ }\def\Jstr{Ann. Sci. \'Ecole Norm. Sup. (4)}%
\def\Vtest{ }\def\Vstr{10}%
\def\Ptest{ }\def\Pstr{441--472}%
\def\Dtest{ }\def\Dstr{1977}%
\def\Astr{\Underlinemark}%
\Refformat\egroup%

\bgroup\Resetstrings%
\def\Ecnt{0}\def\acnt{0}%
\def\Ftest{ }\def\Fstr{6}%
\def\Atest{ }\def\Astr{A\Initper  Borel}%
\def\Ttest{ }\def\Tstr{Linear Algebraic Groups}%
\def\Itest{ }\def\Istr{W.A. Benjamin}%
\def\Ctest{ }\def\Cstr{New York, NY}%
\def\Dtest{ }\def\Dstr{1969}%
\Refformat\egroup%

\bgroup\Resetstrings%
\def\Ecnt{0}\def\acnt{0}%
\def\Ftest{ }\def\Fstr{7}%
\def\Atest{ }\def\Astr{A\Initper  Borel}%
\def\Ttest{ }\def\Tstr{Automorphic $L$\snug-functions}%
\def\Jtest{ }\def\Jstr{Proc. Sympos. Pure Math.}%
\def\Itest{ }\def\Istr{AMS}%
\def\Ctest{ }\def\Cstr{Providence, RI}%
\def\Vtest{ }\def\Vstr{33 part 2}%
\def\Dtest{ }\def\Dstr{1979}%
\def\Ptest{ }\def\Pstr{27--61}%
\def\Astr{\Underlinemark}%
\Refformat\egroup%

\bgroup\Resetstrings%
\def\Ecnt{0}\def\acnt{0}%
\def\Ftest{ }\def\Fstr{8}%
\def\Atest{ }\def\Astr{A\Initper  Borel%
  \Aand H.Jacquet}%
\def\Ttest{ }\def\Tstr{Automorphic forms and automorphic representations}%
\def\Jtest{ }\def\Jstr{Proc. Sympsos. Pure Math.}%
\def\Itest{ }\def\Istr{AMS}%
\def\Vtest{ }\def\Vstr{33 part 1}%
\def\Ctest{ }\def\Cstr{Providence, RI}%
\def\Dtest{ }\def\Dstr{1979}%
\def\Ptest{ }\def\Pstr{189--202}%
\Refformat\egroup%

\bgroup\Resetstrings%
\def\Ecnt{0}\def\acnt{0}%
\def\Ftest{ }\def\Fstr{9}%
\def\Atest{ }\def\Astr{A\Initper  Borel%
  \Aand N\Initper  Wallach}%
\def\Ttest{ }\def\Tstr{Continuous Cohomology, Discrete Subgroups, and
Representations of Reductive Groups}%
\def\Itest{ }\def\Istr{Princeton University Press}%
\def\Stest{ }\def\Sstr{Annals of Math. Studies}%
\def\Ntest{ }\def\Nstr{94}%
\def\Dtest{ }\def\Dstr{1980}%
\def\Ctest{ }\def\Cstr{Princeton, NJ}%
\Refformat\egroup%

\bgroup\Resetstrings%
\def\Ecnt{0}\def\acnt{0}%
\def\Ftest{ }\def\Fstr{10}%
\def\Atest{ }\def\Astr{F\Initper  Bruhat}%
\def\Ttest{ }\def\Tstr{Sur les repr\'esentations induites des groupes de Lie}%
\def\Jtest{ }\def\Jstr{Bull. Soc. Math. France}%
\def\Dtest{ }\def\Dstr{1956}%
\def\Vtest{ }\def\Vstr{84}%
\def\Ptest{ }\def\Pstr{97--205}%
\Refformat\egroup%

\bgroup\Resetstrings%
\def\Ecnt{0}\def\acnt{0}%
\def\Ftest{ }\def\Fstr{11}%
\def\Atest{ }\def\Astr{C\Initper \Initgap J\Initper  Bushnell%
  \Aand P\Initper \Initgap C\Initper  Kutzko}%
\def\Ttest{ }\def\Tstr{To appear}%
\def\Btest{ }\def\Bstr{This volume}%
\Refformat\egroup%

\bgroup\Resetstrings%
\def\Ecnt{0}\def\acnt{0}%
\def\Ftest{ }\def\Fstr{12}%
\def\Atest{ }\def\Astr{C\Initper \Initgap J\Initper  Bushnell%
  \Aand P\Initper \Initgap C\Initper  Kutzko}%
\def\Ttest{ }\def\Tstr{To appear II}%
\def\Btest{ }\def\Bstr{This volume}%
\def\Astr{\Underlinemark}%
\Refformat\egroup%

\bgroup\Resetstrings%
\def\Ecnt{0}\def\acnt{0}%
\def\Ftest{ }\def\Fstr{13}%
\def\Atest{ }\def\Astr{C\Initper \Initgap J\Initper  Bushnell%
  \Aand P\Initper \Initgap C\Initper  Kutzko}%
\def\Ttest{ }\def\Tstr{The Admissible Dual of $GL(N)$ via Open Compact
Subgroups}%
\def\Itest{ }\def\Istr{Princeton Univ. Press}%
\def\Ctest{ }\def\Cstr{Princeton, NJ}%
\def\Stest{ }\def\Sstr{Annals of Mathematics Studies}%
\def\Vtest{ }\def\Vstr{129}%
\def\Dtest{ }\def\Dstr{1993}%
\def\Astr{\Underlinemark}%
\Refformat\egroup%

\bgroup\Resetstrings%
\def\Ecnt{0}\def\acnt{0}%
\def\Ftest{ }\def\Fstr{14}%
\def\Atest{ }\def\Astr{P\Initper  Cartier}%
\def\Ttest{ }\def\Tstr{Representations of $p$\snug-adic groups: A survey}%
\def\Jtest{ }\def\Jstr{Proc. Sympos. Pure Math}%
\def\Itest{ }\def\Istr{AMS}%
\def\Ctest{ }\def\Cstr{Providence, RI}%
\def\Vtest{ }\def\Vstr{33 part 1}%
\def\Ptest{ }\def\Pstr{111--155}%
\def\Dtest{ }\def\Dstr{1979}%
\Refformat\egroup%

\bgroup\Resetstrings%
\def\Ecnt{0}\def\acnt{0}%
\def\Ftest{ }\def\Fstr{15}%
\def\Atest{ }\def\Astr{W\Initper  Casselman}%
\def\Ttest{ }\def\Tstr{Introduction to the theory of admissible representations
of $p$\snug-adic reductive groups}%
\def\Otest{ }\def\Ostr{preprint}%
\Refformat\egroup%

\bgroup\Resetstrings%
\def\Ecnt{0}\def\acnt{0}%
\def\Ftest{ }\def\Fstr{16}%
\def\Atest{ }\def\Astr{W\Initper  Casselman%
  \Aand J\Initper  Shalika}%
\def\Ttest{ }\def\Tstr{The unramified principal series of $p$\snug-adic groups
II,  the Whittaker functions}%
\def\Jtest{ }\def\Jstr{Compositio Math.}%
\def\Vtest{ }\def\Vstr{41}%
\def\Dtest{ }\def\Dstr{1980}%
\def\Ptest{ }\def\Pstr{207--231}%
\Refformat\egroup%

\bgroup\Resetstrings%
\def\Ecnt{0}\def\acnt{0}%
\def\Ftest{ }\def\Fstr{17}%
\def\Atest{ }\def\Astr{S\Initper  Gelbart%
  \Acomma J\Initper  Rogawski%
  \Aandd D\Initper  Soudry}%
\def\Ttest{ }\def\Tstr{On periods of cusp forms and algebraic cycles for
$U(3)$}%
\def\Jtest{ }\def\Jstr{Israel J. Math.}%
\def\Vtest{ }\def\Vstr{83}%
\def\Ptest{ }\def\Pstr{213--252}%
\def\Dtest{ }\def\Dstr{1993}%
\Refformat\egroup%

\bgroup\Resetstrings%
\def\Ecnt{0}\def\acnt{0}%
\def\Ftest{ }\def\Fstr{18}%
\def\Atest{ }\def\Astr{S\Initper \Initgap S\Initper  Gelbart}%
\def\Ttest{ }\def\Tstr{Automorphic Forms on Adele Groups}%
\def\Itest{ }\def\Istr{Princeton University Press}%
\def\Dtest{ }\def\Dstr{1975}%
\def\Ctest{ }\def\Cstr{Princeton, NJ}%
\def\Stest{ }\def\Sstr{Annals of Math. Studies}%
\def\Ntest{ }\def\Nstr{23}%
\Refformat\egroup%

\bgroup\Resetstrings%
\def\Ecnt{0}\def\acnt{0}%
\def\Ftest{ }\def\Fstr{19}%
\def\Atest{ }\def\Astr{S\Initper \Initgap S\Initper  Gelbart}%
\def\Ttest{ }\def\Tstr{An elementary introduction to the Langlands program}%
\def\Jtest{ }\def\Jstr{Bull. Amer. Math. Soc. (N.S.)}%
\def\Vtest{ }\def\Vstr{10}%
\def\Dtest{ }\def\Dstr{1984}%
\def\Ptest{ }\def\Pstr{177--219}%
\def\Astr{\Underlinemark}%
\Refformat\egroup%

\bgroup\Resetstrings%
\def\Ecnt{0}\def\acnt{0}%
\def\Ftest{ }\def\Fstr{20}%
\def\Atest{ }\def\Astr{S\Initper \Initgap S\Initper  Gelbart%
  \Aand A\Initper \Initgap W\Initper  Knapp}%
\def\Ttest{ }\def\Tstr{Irreducible constituents of principal series of
$SL_n(k)$}%
\def\Jtest{ }\def\Jstr{Duke Math. J.}%
\def\Vtest{ }\def\Vstr{48}%
\def\Dtest{ }\def\Dstr{1981}%
\def\Ptest{ }\def\Pstr{313--326}%
\Refformat\egroup%

\bgroup\Resetstrings%
\def\Ecnt{0}\def\acnt{0}%
\def\Ftest{ }\def\Fstr{21}%
\def\Atest{ }\def\Astr{S\Initper \Initgap S\Initper  Gelbart%
  \Aand A\Initper \Initgap W\Initper  Knapp}%
\def\Ttest{ }\def\Tstr{$L$\snug-indistinguishability and $R$ groups for the
special linear group}%
\def\Jtest{ }\def\Jstr{Adv. in Math.}%
\def\Vtest{ }\def\Vstr{43}%
\def\Dtest{ }\def\Dstr{1982}%
\def\Ptest{ }\def\Pstr{101--121}%
\def\Astr{\Underlinemark}%
\Refformat\egroup%

\bgroup\Resetstrings%
\def\Ecnt{0}\def\acnt{0}%
\def\Ftest{ }\def\Fstr{22}%
\def\Atest{ }\def\Astr{S\Initper \Initgap S\Initper  Gelbart%
  \Aand F\Initper  Shahidi}%
\def\Ttest{ }\def\Tstr{Analytic Properties of Automorphic $L$\snug-Functions}%
\def\Itest{ }\def\Istr{Academic Press}%
\def\Ctest{ }\def\Cstr{San Diego, CA}%
\def\Stest{ }\def\Sstr{Perspectives in Mathematics}%
\def\Vtest{ }\def\Vstr{6}%
\def\Dtest{ }\def\Dstr{1988}%
\Refformat\egroup%

\bgroup\Resetstrings%
\def\Ecnt{0}\def\acnt{0}%
\def\Ftest{ }\def\Fstr{23}%
\def\Atest{ }\def\Astr{I\Initper \Initgap M\Initper  Gelfand%
  \Aand D\Initper \Initgap A\Initper  Kazhdan}%
\def\Ttest{ }\def\Tstr{On representations of the group $GL(n, K)$ where $K$ is
a local field}%
\def\Jtest{ }\def\Jstr{Funktsional. Anal. i Prilozhen.}%
\def\Vtest{ }\def\Vstr{6}%
\def\Ntest{ }\def\Nstr{4}%
\def\Dtest{ }\def\Dstr{1972}%
\def\Ptest{ }\def\Pstr{73--74}%
\Refformat\egroup%

\bgroup\Resetstrings%
\def\Ecnt{0}\def\acnt{0}%
\def\Ftest{ }\def\Fstr{24}%
\def\Atest{ }\def\Astr{R\Initper  Godement%
  \Aand H\Initper  Jacquet}%
\def\Ttest{ }\def\Tstr{Zeta Functions of Simple Algebras}%
\def\Itest{ }\def\Istr{Springer-Verlag}%
\def\Ctest{ }\def\Cstr{New York--Heidelberg--Berlin}%
\def\Dtest{ }\def\Dstr{1972}%
\def\Stest{ }\def\Sstr{Lecture Notes in Math.}%
\def\Ntest{ }\def\Nstr{260}%
\Refformat\egroup%

\bgroup\Resetstrings%
\def\Ecnt{0}\def\acnt{0}%
\def\Ftest{ }\def\Fstr{25}%
\def\Atest{ }\def\Astr{D\Initper  Goldberg}%
\def\Ttest{ }\def\Tstr{$R$\snug-groups and elliptic representations for unitary
groups}%
\def\Jtest{ }\def\Jstr{Proc. Amer. Math. Soc.}%
\def\Otest{ }\def\Ostr{to appear}%
\Refformat\egroup%

\bgroup\Resetstrings%
\def\Ecnt{0}\def\acnt{0}%
\def\Ftest{ }\def\Fstr{26}%
\def\Atest{ }\def\Astr{D\Initper  Goldberg}%
\def\Ttest{ }\def\Tstr{Reducibility of induced representations for $Sp(2n)$ and
$SO(n)$}%
\def\Jtest{ }\def\Jstr{Amer. J. Math.}%
\def\Otest{ }\def\Ostr{to appear}%
\def\Astr{\Underlinemark}%
\Refformat\egroup%

\bgroup\Resetstrings%
\def\Ecnt{0}\def\acnt{0}%
\def\Ftest{ }\def\Fstr{27}%
\def\Atest{ }\def\Astr{D\Initper  Goldberg}%
\def\Ttest{ }\def\Tstr{Reducibility of induced representations for classical
$p$\snug-adic groups}%
\def\Dtest{ }\def\Dstr{1991}%
\def\Rtest{ }\def\Rstr{PhD Thesis, University of Maryland}%
\def\Astr{\Underlinemark}%
\Refformat\egroup%

\bgroup\Resetstrings%
\def\Ecnt{0}\def\acnt{0}%
\def\Ftest{ }\def\Fstr{28}%
\def\Atest{ }\def\Astr{D\Initper  Goldberg}%
\def\Ttest{ }\def\Tstr{Reducibility of generalized principal series
representations for $U(2,2)$ via base change}%
\def\Jtest{ }\def\Jstr{Compositio Math.}%
\def\Vtest{ }\def\Vstr{86}%
\def\Ptest{ }\def\Pstr{245--264}%
\def\Dtest{ }\def\Dstr{1993}%
\def\Astr{\Underlinemark}%
\Refformat\egroup%

\bgroup\Resetstrings%
\def\Ecnt{0}\def\acnt{0}%
\def\Ftest{ }\def\Fstr{29}%
\def\Atest{ }\def\Astr{D\Initper  Goldberg}%
\def\Ttest{ }\def\Tstr{$R$\snug-groups and elliptic representations for
$SL_n$}%
\def\Jtest{ }\def\Jstr{Pacific J. Math.}%
\def\Vtest{ }\def\Vstr{165}%
\def\Dtest{ }\def\Dstr{1994}%
\def\Ptest{ }\def\Pstr{77--92}%
\def\Astr{\Underlinemark}%
\Refformat\egroup%

\bgroup\Resetstrings%
\def\Ecnt{0}\def\acnt{0}%
\def\Ftest{ }\def\Fstr{30}%
\def\Atest{ }\def\Astr{D\Initper  Goldberg}%
\def\Ttest{ }\def\Tstr{Some results on reducibility for unitary groups and
local Asai $L$\snug-functions}%
\def\Jtest{ }\def\Jstr{J. Reine Angew. Math.	}%
\def\Vtest{ }\def\Vstr{448}%
\def\Ptest{ }\def\Pstr{65--95}%
\def\Dtest{ }\def\Dstr{1994}%
\def\Astr{\Underlinemark}%
\Refformat\egroup%

\bgroup\Resetstrings%
\def\Ecnt{0}\def\acnt{0}%
\def\Ftest{ }\def\Fstr{31}%
\def\Atest{ }\def\Astr{Harish-Chandra}%
\def\Ttest{ }\def\Tstr{Harmonic analysis on reductive $p$\snug-adic groups}%
\def\Dtest{ }\def\Dstr{1973}%
\def\Jtest{ }\def\Jstr{Proc. Sympos. Pure Math.}%
\def\Itest{ }\def\Istr{AMS}%
\def\Ctest{ }\def\Cstr{Providence, RI}%
\def\Vtest{ }\def\Vstr{26}%
\def\Ptest{ }\def\Pstr{167--192}%
\Refformat\egroup%

\bgroup\Resetstrings%
\def\Ecnt{0}\def\acnt{0}%
\def\Ftest{ }\def\Fstr{32}%
\def\Atest{ }\def\Astr{G\Initper  Henniart}%
\def\Ttest{ }\def\Tstr{To appear}%
\def\Btest{ }\def\Bstr{This volume}%
\Refformat\egroup%

\bgroup\Resetstrings%
\def\Ecnt{0}\def\acnt{0}%
\def\Ftest{ }\def\Fstr{33}%
\def\Atest{ }\def\Astr{R\Initper \Initgap A\Initper  Herb}%
\def\Ttest{ }\def\Tstr{Elliptic representations for $Sp(2n)$ and $SO(n)$}%
\def\Jtest{ }\def\Jstr{Pacific J. Math.}%
\def\Vtest{ }\def\Vstr{161}%
\def\Dtest{ }\def\Dstr{1993}%
\def\Ptest{ }\def\Pstr{347--358}%
\Refformat\egroup%

\bgroup\Resetstrings%
\def\Ecnt{0}\def\acnt{0}%
\def\Ftest{ }\def\Fstr{34}%
\def\Atest{ }\def\Astr{R\Initper  Howe%
  \Aand I\Initper \Initgap I\Initper  Piatetski-Shapiro}%
\def\Ttest{ }\def\Tstr{A counterexample to the "generalized Ramanujan
conjecture" for (quasi-) split groups}%
\def\Jtest{ }\def\Jstr{Proc. Sympsos. Pure Math.}%
\def\Itest{ }\def\Istr{AMS}%
\def\Vtest{ }\def\Vstr{33 part 1}%
\def\Ctest{ }\def\Cstr{Providence, RI}%
\def\Dtest{ }\def\Dstr{1979}%
\def\Ptest{ }\def\Pstr{315--322}%
\Refformat\egroup%

\bgroup\Resetstrings%
\def\Ecnt{0}\def\acnt{0}%
\def\Ftest{ }\def\Fstr{35}%
\def\Atest{ }\def\Astr{J\Initper \Initgap E\Initper  Humphreys}%
\def\Ttest{ }\def\Tstr{Linear Algebraic Groups}%
\def\Itest{ }\def\Istr{Springer-Verlag}%
\def\Ctest{ }\def\Cstr{New York--Heidelberg--Berlin}%
\def\Dtest{ }\def\Dstr{1975}%
\Refformat\egroup%

\bgroup\Resetstrings%
\def\Ecnt{0}\def\acnt{0}%
\def\Ftest{ }\def\Fstr{36}%
\def\Atest{ }\def\Astr{H\Initper  Jacquet}%
\def\Ttest{ }\def\Tstr{Repr\'esentations des groupes lin\'eares
$p$\snug-adiques}%
\def\Btest{ }\def\Bstr{Theory of Group representations and Fourier Analysis}%
\def\Itest{ }\def\Istr{CIME}%
\def\Ctest{ }\def\Cstr{Montecatini}%
\def\Dtest{ }\def\Dstr{1971}%
\def\Ptest{ }\def\Pstr{121--220}%
\Refformat\egroup%

\bgroup\Resetstrings%
\def\Ecnt{0}\def\acnt{0}%
\def\Ftest{ }\def\Fstr{37}%
\def\Atest{ }\def\Astr{H\Initper  Jacquet}%
\def\Ttest{ }\def\Tstr{Generic representations}%
\def\Btest{ }\def\Bstr{Non Commutatuve Harmonic Analysis}%
\def\Itest{ }\def\Istr{Springer-Verlag}%
\def\Ctest{ }\def\Cstr{New York--Heidelberg--Berlin}%
\def\Stest{ }\def\Sstr{Lecture Notes in Mathematics}%
\def\Ntest{ }\def\Nstr{587}%
\def\Dtest{ }\def\Dstr{1977}%
\def\Ptest{ }\def\Pstr{91--101}%
\def\Astr{\Underlinemark}%
\Refformat\egroup%

\bgroup\Resetstrings%
\def\Ecnt{0}\def\acnt{0}%
\def\Ftest{ }\def\Fstr{38}%
\def\Atest{ }\def\Astr{C\Initper \Initgap D\Initper  Keys}%
\def\Ttest{ }\def\Tstr{On the decomposition of reducible principal series
representations of $p$\snug-adic Chevalley groups}%
\def\Dtest{ }\def\Dstr{1982}%
\def\Jtest{ }\def\Jstr{Pacific J. Math.}%
\def\Vtest{ }\def\Vstr{101}%
\def\Ptest{ }\def\Pstr{351--388}%
\Refformat\egroup%

\bgroup\Resetstrings%
\def\Ecnt{0}\def\acnt{0}%
\def\Ftest{ }\def\Fstr{39}%
\def\Atest{ }\def\Astr{C\Initper \Initgap D\Initper  Keys}%
\def\Ttest{ }\def\Tstr{L-indistinguishability and R-groups for quasi split
groups: unitary groups in even dimension}%
\def\Jtest{ }\def\Jstr{Ann. Sci. \'Ecole Norm. Sup. (4)}%
\def\Dtest{ }\def\Dstr{1987}%
\def\Vtest{ }\def\Vstr{20}%
\def\Ptest{ }\def\Pstr{31--64}%
\def\Astr{\Underlinemark}%
\Refformat\egroup%

\bgroup\Resetstrings%
\def\Ecnt{0}\def\acnt{0}%
\def\Ftest{ }\def\Fstr{40}%
\def\Atest{ }\def\Astr{A\Initper \Initgap W\Initper  Knapp%
  \Aand E\Initper \Initgap M\Initper  Stein}%
\def\Ttest{ }\def\Tstr{Intertwining operators for semisimplelie groups}%
\def\Jtest{ }\def\Jstr{Ann. of Math. (2)}%
\def\Dtest{ }\def\Dstr{1971}%
\def\Vtest{ }\def\Vstr{93}%
\def\Ptest{ }\def\Pstr{489--578}%
\Refformat\egroup%

\bgroup\Resetstrings%
\def\Ecnt{0}\def\acnt{0}%
\def\Ftest{ }\def\Fstr{41}%
\def\Atest{ }\def\Astr{A\Initper \Initgap W\Initper  Knapp%
  \Aand E\Initper \Initgap M\Initper  Stein}%
\def\Ttest{ }\def\Tstr{Irreducibility theorems for the principal series}%
\def\Btest{ }\def\Bstr{Conference on Harmonic Analysis}%
\def\Itest{ }\def\Istr{Springer-Verlag}%
\def\Ctest{ }\def\Cstr{New York--Heidelberg--Berlin}%
\def\Stest{ }\def\Sstr{Lecture Notes in Mathematics}%
\def\Dtest{ }\def\Dstr{1972}%
\def\Ntest{ }\def\Nstr{266}%
\def\Ptest{ }\def\Pstr{197--214}%
\def\Astr{\Underlinemark}%
\Refformat\egroup%

\bgroup\Resetstrings%
\def\Ecnt{0}\def\acnt{0}%
\def\Ftest{ }\def\Fstr{42}%
\def\Atest{ }\def\Astr{A\Initper \Initgap W\Initper  Knapp%
  \Aand G\Initper  Zuckerman}%
\def\Ttest{ }\def\Tstr{Classification of irreducible tempered representations
of semisimple Lie groups}%
\def\Jtest{ }\def\Jstr{Proc. Nat. Acad. Sci. U.S.A.}%
\def\Vtest{ }\def\Vstr{73}%
\def\Ntest{ }\def\Nstr{7}%
\def\Ptest{ }\def\Pstr{2178--2180}%
\def\Dtest{ }\def\Dstr{1976}%
\Refformat\egroup%

\bgroup\Resetstrings%
\def\Ecnt{0}\def\acnt{0}%
\def\Ftest{ }\def\Fstr{43}%
\def\Atest{ }\def\Astr{A\Initper \Initgap W\Initper  Knapp%
  \Aand G\Initper  Zuckerman}%
\def\Ttest{ }\def\Tstr{Multiplicity one fails for $p$\snug-adic unitary
principal series}%
\def\Jtest{ }\def\Jstr{Hiroshima Math. J.}%
\def\Vtest{ }\def\Vstr{10}%
\def\Dtest{ }\def\Dstr{1980}%
\def\Ptest{ }\def\Pstr{295--309}%
\def\Astr{\Underlinemark}%
\Refformat\egroup%

\bgroup\Resetstrings%
\def\Ecnt{0}\def\acnt{0}%
\def\Ftest{ }\def\Fstr{44}%
\def\Atest{ }\def\Astr{R\Initper  Kottwitz%
  \Aand D\Initper  Shelstad}%
\def\Ttest{ }\def\Tstr{Twisted Endoscopy II: Basic global theory}%
\def\Otest{ }\def\Ostr{preprint}%
\Refformat\egroup%

\bgroup\Resetstrings%
\def\Ecnt{0}\def\acnt{0}%
\def\Ftest{ }\def\Fstr{45}%
\def\Atest{ }\def\Astr{R\Initper  Kottwitz%
  \Aand D\Initper  Shelstad}%
\def\Ttest{ }\def\Tstr{Twisted Endoscopy I: Definitions, norm mappings and
transfer factors}%
\def\Otest{ }\def\Ostr{preprint}%
\def\Astr{\Underlinemark}%
\Refformat\egroup%

\bgroup\Resetstrings%
\def\Ecnt{0}\def\acnt{0}%
\def\Ftest{ }\def\Fstr{46}%
\def\Atest{ }\def\Astr{R\Initper \Initgap A\Initper  Kunze%
  \Aand E\Initper \Initgap M\Initper  Stein}%
\def\Ttest{ }\def\Tstr{Uniformly bounded representations and harmonic analysis
of the $2\times 2$ real unimodular group}%
\def\Jtest{ }\def\Jstr{Amer. J. Math.}%
\def\Vtest{ }\def\Vstr{82}%
\def\Ptest{ }\def\Pstr{1--62}%
\def\Dtest{ }\def\Dstr{1960}%
\Refformat\egroup%

\bgroup\Resetstrings%
\def\Ecnt{0}\def\acnt{0}%
\def\Ftest{ }\def\Fstr{47}%
\def\Atest{ }\def\Astr{R\Initper \Initgap A\Initper  Kunze%
  \Aand E\Initper \Initgap M\Initper  Stein}%
\def\Ttest{ }\def\Tstr{Uniformly bounded representations II, analytic
continuation of the principal series representations of the $n\times n$ complex
unimodular group}%
\def\Jtest{ }\def\Jstr{Amer. J. Math.}%
\def\Vtest{ }\def\Vstr{83}%
\def\Ptest{ }\def\Pstr{723--786}%
\def\Dtest{ }\def\Dstr{1961}%
\def\Astr{\Underlinemark}%
\Refformat\egroup%

\bgroup\Resetstrings%
\def\Ecnt{0}\def\acnt{0}%
\def\Ftest{ }\def\Fstr{48}%
\def\Atest{ }\def\Astr{R\Initper \Initgap A\Initper  Kunze%
  \Aand E\Initper \Initgap M\Initper  Stein}%
\def\Ttest{ }\def\Tstr{Uniformly bounded representations III, intertwining
operators for the principal series on semisimple groups}%
\def\Jtest{ }\def\Jstr{Amer. J. Math.}%
\def\Vtest{ }\def\Vstr{89}%
\def\Ptest{ }\def\Pstr{385--442}%
\def\Dtest{ }\def\Dstr{1967}%
\def\Astr{\Underlinemark}%
\Refformat\egroup%

\bgroup\Resetstrings%
\def\Ecnt{0}\def\acnt{0}%
\def\Ftest{ }\def\Fstr{49}%
\def\Atest{ }\def\Astr{R\Initper \Initgap A\Initper  Kunze%
  \Aand E\Initper \Initgap M\Initper  Stein}%
\def\Ttest{ }\def\Tstr{Uniformly bounded representations IV, analytic
continuation of the principal series for complex classical groups of types
$B_n,\ C_n,\ D_n$}%
\def\Jtest{ }\def\Jstr{Adv. in Math.}%
\def\Vtest{ }\def\Vstr{11}%
\def\Ptest{ }\def\Pstr{1--71}%
\def\Dtest{ }\def\Dstr{1973}%
\def\Astr{\Underlinemark}%
\Refformat\egroup%

\bgroup\Resetstrings%
\def\Ecnt{0}\def\acnt{0}%
\def\Ftest{ }\def\Fstr{50}%
\def\Atest{ }\def\Astr{J\Initper \Initgap P\Initper  Labesse%
  \Aand R\Initper \Initgap P\Initper  Langlands}%
\def\Ttest{ }\def\Tstr{$L$\snug-indistinguishability for $SL(2)$}%
\def\Jtest{ }\def\Jstr{Canad. J. Math.}%
\def\Vtest{ }\def\Vstr{31}%
\def\Dtest{ }\def\Dstr{1979}%
\def\Ptest{ }\def\Pstr{726--785}%
\Refformat\egroup%

\bgroup\Resetstrings%
\def\Ecnt{0}\def\acnt{0}%
\def\Ftest{ }\def\Fstr{51}%
\def\Atest{ }\def\Astr{R\Initper \Initgap P\Initper  Langlands}%
\def\Ttest{ }\def\Tstr{Problems in the theory of automorphic froms}%
\def\Btest{ }\def\Bstr{Lecture notes in Mathematics}%
\def\Vtest{ }\def\Vstr{170}%
\def\Dtest{ }\def\Dstr{1970}%
\def\Itest{ }\def\Istr{Springer-Verlag}%
\def\Ctest{ }\def\Cstr{New York--Heidelberg--Berlin}%
\def\Ptest{ }\def\Pstr{18--86}%
\Refformat\egroup%

\bgroup\Resetstrings%
\def\Ecnt{0}\def\acnt{0}%
\def\Ftest{ }\def\Fstr{52}%
\def\Atest{ }\def\Astr{R\Initper \Initgap P\Initper  Langlands}%
\def\Ttest{ }\def\Tstr{Euler Products}%
\def\Otest{ }\def\Ostr{Yale University}%
\def\Dtest{ }\def\Dstr{1971}%
\def\Astr{\Underlinemark}%
\Refformat\egroup%

\bgroup\Resetstrings%
\def\Ecnt{0}\def\acnt{0}%
\def\Ftest{ }\def\Fstr{53}%
\def\Atest{ }\def\Astr{R\Initper \Initgap P\Initper  Langlands}%
\def\Ttest{ }\def\Tstr{On the functional equation satisfied by Eisenstein
series}%
\def\Itest{ }\def\Istr{Springer-Verlag}%
\def\Ctest{ }\def\Cstr{New York--Heidelberg--Berlin}%
\def\Stest{ }\def\Sstr{Lecture Notes in Mathematics }%
\def\Ntest{ }\def\Nstr{544}%
\def\Dtest{ }\def\Dstr{1976}%
\def\Astr{\Underlinemark}%
\Refformat\egroup%

\bgroup\Resetstrings%
\def\Ecnt{0}\def\acnt{0}%
\def\Ftest{ }\def\Fstr{54}%
\def\Atest{ }\def\Astr{R\Initper \Initgap P\Initper  Langlands}%
\def\Ttest{ }\def\Tstr{On the notion of an automorphic representation}%
\def\Jtest{ }\def\Jstr{Proc. Sympsos. Pure Math.}%
\def\Itest{ }\def\Istr{AMS}%
\def\Vtest{ }\def\Vstr{33 part 1}%
\def\Ctest{ }\def\Cstr{Providence, RI}%
\def\Dtest{ }\def\Dstr{1979}%
\def\Ptest{ }\def\Pstr{203--207}%
\def\Astr{\Underlinemark}%
\Refformat\egroup%

\bgroup\Resetstrings%
\def\Ecnt{0}\def\acnt{0}%
\def\Ftest{ }\def\Fstr{55}%
\def\Atest{ }\def\Astr{R\Initper \Initgap P\Initper  Langlands%
  \Aand D\Initper  Shelstad}%
\def\Ttest{ }\def\Tstr{On the definition of transfer factors}%
\def\Jtest{ }\def\Jstr{Math. Ann.}%
\def\Vtest{ }\def\Vstr{278}%
\def\Dtest{ }\def\Dstr{1987}%
\def\Ptest{ }\def\Pstr{219--271}%
\Refformat\egroup%

\bgroup\Resetstrings%
\def\Ecnt{0}\def\acnt{0}%
\def\Ftest{ }\def\Fstr{56}%
\def\Atest{ }\def\Astr{L\Initper  Morris}%
\def\Ttest{ }\def\Tstr{Reductive groups over local fields}%
\def\Btest{ }\def\Bstr{This Volume}%
\def\Rtest{ }\def\Rstr{To appear}%
\Refformat\egroup%

\bgroup\Resetstrings%
\def\Ecnt{0}\def\acnt{0}%
\def\Ftest{ }\def\Fstr{57}%
\def\Atest{ }\def\Astr{L\Initper  Morris}%
\def\Ttest{ }\def\Tstr{Reductive groups}%
\def\Btest{ }\def\Bstr{This volume}%
\def\Rtest{ }\def\Rstr{To appear}%
\def\Astr{\Underlinemark}%
\Refformat\egroup%

\bgroup\Resetstrings%
\def\Ecnt{0}\def\acnt{0}%
\def\Ftest{ }\def\Fstr{58}%
\def\Atest{ }\def\Astr{G\Initper \Initgap I\Initper  Ol'{\v{s}}anski{\v{i}}}%
\def\Ttest{ }\def\Tstr{Intertwining operators and complementary series in the
class of representations induced from parabolic subgroups of the general linear
group over a locally compact division algebra}%
\def\Jtest{ }\def\Jstr{Math. USSR-Sb.}%
\def\Dtest{ }\def\Dstr{1974}%
\def\Vtest{ }\def\Vstr{22}%
\def\Ptest{ }\def\Pstr{217--254}%
\Refformat\egroup%

\bgroup\Resetstrings%
\def\Ecnt{0}\def\acnt{0}%
\def\Ftest{ }\def\Fstr{59}%
\def\Atest{ }\def\Astr{F\Initper  Rodier}%
\def\Ttest{ }\def\Tstr{Whittaker models for admissible representations of
reductive $p$\snug-adic split groups}%
\def\Jtest{ }\def\Jstr{Proc. Sympos. Pure Math.}%
\def\Itest{ }\def\Istr{AMS}%
\def\Ctest{ }\def\Cstr{Providence, RI}%
\def\Vtest{ }\def\Vstr{26}%
\def\Ptest{ }\def\Pstr{425--430}%
\def\Dtest{ }\def\Dstr{1973}%
\Refformat\egroup%

\bgroup\Resetstrings%
\def\Ecnt{0}\def\acnt{0}%
\def\Ftest{ }\def\Fstr{60}%
\def\Atest{ }\def\Astr{J\Initper \Initgap D\Initper  Rogawski}%
\def\Ttest{ }\def\Tstr{Automorphic Representations of Unitary Groups in Three
Variables}%
\def\Itest{ }\def\Istr{Princeton University Press}%
\def\Stest{ }\def\Sstr{Annals of Math. Studies}%
\def\Ntest{ }\def\Nstr{123}%
\def\Dtest{ }\def\Dstr{1990}%
\def\Ctest{ }\def\Cstr{Princeton, NJ}%
\Refformat\egroup%

\bgroup\Resetstrings%
\def\Ecnt{0}\def\acnt{0}%
\def\Ftest{ }\def\Fstr{61}%
\def\Atest{ }\def\Astr{P\Initper \Initgap J\Initper  Sally-Jr.}%
\def\Ttest{ }\def\Tstr{Unitary and uniformly bounded representations of the two
by two unimodular group over local fields}%
\def\Jtest{ }\def\Jstr{Amer. J. Math.}%
\def\Dtest{ }\def\Dstr{1968}%
\def\Vtest{ }\def\Vstr{90}%
\def\Ptest{ }\def\Pstr{406--443}%
\Refformat\egroup%

\bgroup\Resetstrings%
\def\Ecnt{0}\def\acnt{0}%
\def\Ftest{ }\def\Fstr{62}%
\def\Atest{ }\def\Astr{G\Initper  Schiffman}%
\def\Ttest{ }\def\Tstr{Int\'egrales d'entrelacement et functions de Whittaker}%
\def\Jtest{ }\def\Jstr{Bull. Soc. Math. France}%
\def\Vtest{ }\def\Vstr{99}%
\def\Dtest{ }\def\Dstr{1971}%
\def\Ptest{ }\def\Pstr{3--72}%
\Refformat\egroup%

\bgroup\Resetstrings%
\def\Ecnt{0}\def\acnt{0}%
\def\Ftest{ }\def\Fstr{63}%
\def\Atest{ }\def\Astr{F\Initper  Shahidi}%
\def\Ttest{ }\def\Tstr{The notion of norm and the representation theory of
orthogonal groups}%
\def\Jtest{ }\def\Jstr{Invent. Math.}%
\def\Otest{ }\def\Ostr{To appear}%
\Refformat\egroup%

\bgroup\Resetstrings%
\def\Ecnt{0}\def\acnt{0}%
\def\Ftest{ }\def\Fstr{64}%
\def\Atest{ }\def\Astr{F\Initper  Shahidi}%
\def\Ttest{ }\def\Tstr{On certain $L$\snug-functions}%
\def\Jtest{ }\def\Jstr{Amer. J. Math.}%
\def\Vtest{ }\def\Vstr{103}%
\def\Ntest{ }\def\Nstr{2}%
\def\Ptest{ }\def\Pstr{297--355}%
\def\Dtest{ }\def\Dstr{1981}%
\def\Astr{\Underlinemark}%
\Refformat\egroup%

\bgroup\Resetstrings%
\def\Ecnt{0}\def\acnt{0}%
\def\Ftest{ }\def\Fstr{65}%
\def\Atest{ }\def\Astr{F\Initper  Shahidi}%
\def\Ttest{ }\def\Tstr{Fourier transforms of intertwining operators and
Placherel measures for  $GL(n)$}%
\def\Jtest{ }\def\Jstr{Amer. J. Math.}%
\def\Vtest{ }\def\Vstr{106}%
\def\Dtest{ }\def\Dstr{1984}%
\def\Ptest{ }\def\Pstr{67--111}%
\def\Astr{\Underlinemark}%
\Refformat\egroup%

\bgroup\Resetstrings%
\def\Ecnt{0}\def\acnt{0}%
\def\Ftest{ }\def\Fstr{66}%
\def\Atest{ }\def\Astr{F\Initper  Shahidi}%
\def\Ttest{ }\def\Tstr{On the Ramanujan conjecture and finiteness of poles for
certain $L$\snug-functions}%
\def\Jtest{ }\def\Jstr{Ann. of Math. (2)}%
\def\Vtest{ }\def\Vstr{127}%
\def\Dtest{ }\def\Dstr{1988}%
\def\Ptest{ }\def\Pstr{547--584}%
\def\Astr{\Underlinemark}%
\Refformat\egroup%

\bgroup\Resetstrings%
\def\Ecnt{0}\def\acnt{0}%
\def\Ftest{ }\def\Fstr{67}%
\def\Atest{ }\def\Astr{F\Initper  Shahidi}%
\def\Ttest{ }\def\Tstr{A proof of Langlands conjecture for Plancherel measures;
complementary series for $p$\snug-adic groups}%
\def\Jtest{ }\def\Jstr{Ann. of Math. (2)}%
\def\Dtest{ }\def\Dstr{1990}%
\def\Vtest{ }\def\Vstr{132}%
\def\Ptest{ }\def\Pstr{273--330}%
\def\Astr{\Underlinemark}%
\Refformat\egroup%

\bgroup\Resetstrings%
\def\Ecnt{0}\def\acnt{0}%
\def\Ftest{ }\def\Fstr{68}%
\def\Atest{ }\def\Astr{F\Initper  Shahidi}%
\def\Ttest{ }\def\Tstr{Twisted endoscopy and reducibility of induced
representations for $p$\snug-adic groups}%
\def\Jtest{ }\def\Jstr{Duke Math. J.}%
\def\Vtest{ }\def\Vstr{66}%
\def\Dtest{ }\def\Dstr{1992}%
\def\Ptest{ }\def\Pstr{1--41}%
\def\Astr{\Underlinemark}%
\Refformat\egroup%

\bgroup\Resetstrings%
\def\Ecnt{2}\def\acnt{0}%
\def\Ftest{ }\def\Fstr{69}%
\def\Atest{ }\def\Astr{F\Initper  Shahidi}%
\def\Ttest{ }\def\Tstr{$L$\snug-functions and representation theory of
$p$\snug-adic groups}%
\def\Btest{ }\def\Bstr{$p$\snug-adic Methods and Their Applications}%
\def\Etest{ }\def\Estr{A\Initper \Initgap J\Initper  Baker%
  \Eand R\Initper \Initgap J\Initper  Plymen}%
\def\Itest{ }\def\Istr{Clarendon Press}%
\def\Ctest{ }\def\Cstr{Oxford New York Tokyo}%
\def\Dtest{ }\def\Dstr{1992}%
\def\Ptest{ }\def\Pstr{91--112}%
\def\Astr{\Underlinemark}%
\Refformat\egroup%

\bgroup\Resetstrings%
\def\Ecnt{0}\def\acnt{0}%
\def\Ftest{ }\def\Fstr{70}%
\def\Atest{ }\def\Astr{J\Initper  Shalika}%
\def\Ttest{ }\def\Tstr{The multiplicity one theorem for $GL_n$}%
\def\Jtest{ }\def\Jstr{Ann. of Math.}%
\def\Vtest{ }\def\Vstr{100}%
\def\Dtest{ }\def\Dstr{1974}%
\def\Ptest{ }\def\Pstr{171--193}%
\Refformat\egroup%

\bgroup\Resetstrings%
\def\Ecnt{0}\def\acnt{0}%
\def\Ftest{ }\def\Fstr{71}%
\def\Atest{ }\def\Astr{A\Initper \Initgap J\Initper  Silberger}%
\def\Ttest{ }\def\Tstr{The Knapp-Stein dimension theorem for $p$\snug-adic
groups}%
\def\Jtest{ }\def\Jstr{Proc. Amer. Math. Soc.}%
\def\Vtest{ }\def\Vstr{68}%
\def\Dtest{ }\def\Dstr{1978}%
\def\Ptest{ }\def\Pstr{243--246}%
\Refformat\egroup%

\bgroup\Resetstrings%
\def\Ecnt{0}\def\acnt{0}%
\def\Ftest{ }\def\Fstr{72}%
\def\Atest{ }\def\Astr{A\Initper \Initgap J\Initper  Silberger}%
\def\Ttest{ }\def\Tstr{Introduction to Harmonic Analysis on Reductive
$p$\snug-adic Groups}%
\def\Itest{ }\def\Istr{Princeton University Press}%
\def\Ctest{ }\def\Cstr{Princeton, NJ}%
\def\Stest{ }\def\Sstr{Mathematical Notes}%
\def\Ntest{ }\def\Nstr{23}%
\def\Dtest{ }\def\Dstr{1979}%
\def\Astr{\Underlinemark}%
\Refformat\egroup%

\bgroup\Resetstrings%
\def\Ecnt{0}\def\acnt{0}%
\def\Ftest{ }\def\Fstr{73}%
\def\Atest{ }\def\Astr{A\Initper \Initgap J\Initper  Silberger}%
\def\Ttest{ }\def\Tstr{The Knapp-Stein dimension theorem for $p$\snug-adic
groups. Correction}%
\def\Jtest{ }\def\Jstr{Proc. Amer. Math. Soc.}%
\def\Vtest{ }\def\Vstr{76}%
\def\Dtest{ }\def\Dstr{1979}%
\def\Ptest{ }\def\Pstr{169--170}%
\def\Astr{\Underlinemark}%
\Refformat\egroup%

\bgroup\Resetstrings%
\def\Ecnt{0}\def\acnt{0}%
\def\Ftest{ }\def\Fstr{74}%
\def\Atest{ }\def\Astr{R\Initper  Steinberg}%
\def\Ttest{ }\def\Tstr{Lectures on Chevalley Groups}%
\def\Jtest{ }\def\Jstr{Yale University Lecture Notes}%
\def\Dtest{ }\def\Dstr{1967}%
\Refformat\egroup%

\bgroup\Resetstrings%
\def\Ecnt{0}\def\acnt{0}%
\def\Ftest{ }\def\Fstr{75}%
\def\Atest{ }\def\Astr{M\Initper  Tadic}%
\def\Ttest{ }\def\Tstr{Representations of $p$\snug-adic symplectic groups}%
\def\Jtest{ }\def\Jstr{Compositio Math.}%
\def\Otest{ }\def\Ostr{to appear}%
\Refformat\egroup%

\bgroup\Resetstrings%
\def\Ecnt{2}\def\acnt{0}%
\def\Ftest{ }\def\Fstr{76}%
\def\Atest{ }\def\Astr{M\Initper  Tadic}%
\def\Ttest{ }\def\Tstr{On Jacquet modules of induced representations of
$p$\snug-adic symplectic groups}%
\def\Btest{ }\def\Bstr{Harmonic Analysis on Reducitve Groups}%
\def\Dtest{ }\def\Dstr{1991}%
\def\Itest{ }\def\Istr{Birkh{\"a}user Boston}%
\def\Ctest{ }\def\Cstr{Cambridge, MA}%
\def\Ptest{ }\def\Pstr{57--78}%
\def\Etest{ }\def\Estr{W\Initper  Barker%
  \Eand P\Initper  Sally}%
\def\Astr{\Underlinemark}%
\Refformat\egroup%

\bgroup\Resetstrings%
\def\Ecnt{0}\def\acnt{0}%
\def\Ftest{ }\def\Fstr{77}%
\def\Atest{ }\def\Astr{M\Initper  Tadic}%
\def\Ttest{ }\def\Tstr{Structure arising from induction and Jacquet modules of
representations of classical $p$\snug-adic groups}%
\def\Jtest{ }\def\Jstr{Mathematica Gottengensis}%
\def\Vtest{ }\def\Vstr{41}%
\def\Ntest{ }\def\Nstr{9}%
\def\Dtest{ }\def\Dstr{1992}%
\def\Astr{\Underlinemark}%
\Refformat\egroup%

\bgroup\Resetstrings%
\def\Ecnt{0}\def\acnt{0}%
\def\Ftest{ }\def\Fstr{78}%
\def\Atest{ }\def\Astr{M\Initper  Tadic}%
\def\Ttest{ }\def\Tstr{Notes on representations of non-archimedean $SL(n)$}%
\def\Jtest{ }\def\Jstr{Pacific J. Math.}%
\def\Vtest{ }\def\Vstr{152}%
\def\Dtest{ }\def\Dstr{1992}%
\def\Ptest{ }\def\Pstr{375--396}%
\def\Astr{\Underlinemark}%
\Refformat\egroup%

\bgroup\Resetstrings%
\def\Ecnt{2}\def\acnt{0}%
\def\Ftest{ }\def\Fstr{79}%
\def\Atest{ }\def\Astr{J\Initper  Tate}%
\def\Ttest{ }\def\Tstr{Fourier analysis in number fields and Hecke's zeta
functions}%
\def\Btest{ }\def\Bstr{Algebraic Number Theory}%
\def\Etest{ }\def\Estr{J\Initper \Initgap W\Initper  Cassels%
  \Eand A\Initper  Fr\"olich}%
\def\Itest{ }\def\Istr{Academic Press}%
\def\Dtest{ }\def\Dstr{1968}%
\def\Ptest{ }\def\Pstr{305--347}%
\Refformat\egroup%

\bgroup\Resetstrings%
\def\Ecnt{0}\def\acnt{0}%
\def\Ftest{ }\def\Fstr{80}%
\def\Atest{ }\def\Astr{J\Initper  Tate}%
\def\Ttest{ }\def\Tstr{Number theoretic background}%
\def\Jtest{ }\def\Jstr{Proc. Sympos. Pure Math.}%
\def\Itest{ }\def\Istr{AMS}%
\def\Ctest{ }\def\Cstr{Providence, RI}%
\def\Vtest{ }\def\Vstr{33 part 2}%
\def\Dtest{ }\def\Dstr{1979}%
\def\Ptest{ }\def\Pstr{3--26}%
\def\Astr{\Underlinemark}%
\Refformat\egroup%

\bgroup\Resetstrings%
\def\Ecnt{0}\def\acnt{0}%
\def\Ftest{ }\def\Fstr{81}%
\def\Atest{ }\def\Astr{N\Initper  Winarsky}%
\def\Ttest{ }\def\Tstr{Reducibility of principal series representations of
$p$\snug-adic Chevalley groups}%
\def\Jtest{ }\def\Jstr{Amer. J. Math.}%
\def\Vtest{ }\def\Vstr{100}%
\def\Dtest{ }\def\Dstr{1978}%
\def\Ntest{ }\def\Nstr{5}%
\def\Ptest{ }\def\Pstr{941--956}%
\Refformat\egroup%

\bgroup\Resetstrings%
\def\Ecnt{0}\def\acnt{0}%
\def\Ftest{ }\def\Fstr{82}%
\def\Atest{ }\def\Astr{A\Initper \Initgap V\Initper  Zelevinsky}%
\def\Ttest{ }\def\Tstr{Induced representations of reductive $p$\snug-adic
groups II, on irreducible representations of $GL(n)$}%
\def\Jtest{ }\def\Jstr{Ann. Sci. \'Ecole Norm. Sup. (4)}%
\def\Vtest{ }\def\Vstr{13}%
\def\Dtest{ }\def\Dstr{1980}%
\def\Ptest{ }\def\Pstr{165--210}%
\Refformat\egroup%

\endRefs
\enddocument